\documentclass[12pt]{amsart} 
\usepackage{amsthm}  
\usepackage{epsf}  
\usepackage{amssymb}  
\usepackage{latexsym}  
\usepackage{amsmath}  
\newtheorem{theorem}{Theorem}[section] 
\newtheorem{thm}[theorem]{Theorem}
  
\newtheorem{conj}[theorem]{Conjecture}  
\newtheorem{cor}[theorem]{Corollary}

\newtheorem{prop}[theorem]{Proposition}  
\newtheorem{rem}[theorem]{Remark}  

\def\R{{\Bbb R}}  
  
\def\C{{\Bbb C}}  
  
\def\Z{{\Bbb Z}}  
\def\P{{\Bbb P}}  
\def\F{{\Bbb F}}  
  
\def\cpt{\Bbb C \Bbb P^2}

\def \ta{\tau}  
\def \ta1{\tau_1}

\def \De{\Delta}  
\def \de{\delta}  
\def \eps{\epsilon}

\def \la{\lambda}

\def\pcpt{\pi_1(\cpt - S, *)}

\def \zovera {  
    \mathop{\lower 10pt \hbox{${\buildrel{\displaystyle\bar{z}} \over {\scriptstyle{(a)}}} $}}  
    {\lower 4pt \hbox{${\scriptstyle{ij}}$}} 
} 

\def\pcpt{\pi_1 (\C\P^2 - S,*)}

\long\def\forget#1\forgotten{}
  
\newcommand\begintable[1][] {{}}

\newif\ifXY 
\XYfalse    
\newif\ifbigmatrices  
\bigmatricestrue 
  
 \ifXY  
 \usepackage{xy}  
 \fi  
 \ifXY  
 \xyoption{all}  
 \fi  
  
\begin{document}  
  
\title[$\pi_1$ of arrangements with singularities up to order 6] 
{Fundamental groups of Tangent Conic-Line arrangements 
with Singularities up to order 6}
  
\author[Meirav Amram, David Garber and Mina Teicher]{Meirav Amram$^1$, David Garber$^2$ and Mina Teicher}  
  
\stepcounter{footnote} 
\footnotetext{Partially supported by the Minerva Foundation of Germany, the DAAD 
fellowship (Germany), EAGER (EU Network, HPRN-CT-2009-00099) and 
the Golda Meir postdoctoral fellowship.}
\stepcounter{footnote} 
\footnotetext{Partially supported by Lady Davis and Golda Meir postdoctoral 
fellowships.} 

\address{Meirav Amram, Einstein Institute of Mathematics, 
The Hebrew University, Givat Ram, 91904 Jerusalem, Israel}
\email{ameirav@math.huji.ac.il}

\address{David Garber, Einstein Institute of Mathematics, 
The Hebrew University, Givat Ram, 91904 Jerusalem, Israel, and 
Department of Applied Mathematics, Holon Academic Institute of Technology, 
52 Golomb St., PO Box 305, 58102 Holon, Israel}

\email{garber@math.huji.ac.il,garber@hait.ac.il}

\address{Mina Teicher, Department of Mathematics, Bar-Ilan University, 52900
Ramat Gan, Israel}  
  
\email{teicher@macs.biu.ac.il}

  
 \renewcommand{\subjclassname}{%
       \textup{2000} Mathematics Subject Classification}  
  
  
 \date{\today}

\begin{abstract}  
We list all the possible fundamental groups of the complements of real conic-line 
arrangements with two conics which are tangent to each other at two points, 
with up to two additional lines. 

For the computations we use the topological local braid monodromies  
and the techniques of Moishezon-Teicher and van-Kampen.

We also include some conjectures concerning the connection  
between the presentation of the fundamental group of 
the complements and the geometry of an interesting family of conic line 
arrangements.  
\end{abstract}  
  
\maketitle  
  
\section{Introduction}\label{intro}  
  
In this paper we compute and list the fundamental groups of complements of
all real conic-line arrangements in $\cpt$ with two conics, which are tangent to each other 
at two points, with up to two additional lines in any position.

Algorithmically, this paper uses the local computations 
(local braid monodromies and their induced relations), 
the braid monodromy techniques of Moishezon-Teicher
(see \cite{GaTe}, \cite{Mo1}, \cite{MoTe1}, \cite{MoTe2}, 
\cite{MoTe3} and \cite{MoTe4}),  
the Enriques - van Kampen Theorem (see \cite{vK})
and some group calculations for studying the fundamental groups.
See \cite{AmTe} for detailed exposition of these techniques. 

These arrangements may appear in a branch curve of a generic projection 
to $\cpt$ of a surface of general type (see for example \cite{Hi}). 
 
The main results of this paper are as follows 
(where $e$ is the unit element in the group and $\F_2$ is the free group with two generators):  
\begin{prop}\label{th1-2}  
Let $S$ be a curve in $\cpt$ composed of two tangent conics. Then:  
$$\pcpt \cong \langle x_1,x_2\ |\ (x_1 x_2)^2=(x_2 x_1)^2=e \rangle$$
\end{prop}  
  
\begin{prop}\label{th1}
There are three possible fundamental groups for 
a conic-line arrangement which consists of two tangent conics and an additional line:
\begin{enumerate}
\item $\Z \oplus \langle x_1,x_2 \ | \ (x_1 x_2)^2=(x_2 x_1)^2=e \rangle$.
\item $\langle x_1,x_2 \ | \ (x_1 x_2)^2 = (x_2 x_1)^2 \rangle$.
\item $\F_2$
\end{enumerate}
\end{prop}  

\begin{prop}\label{th2}
There are five possible fundamental groups for 
a conic-line arrangement which consists of two tangent conics and two additional lines:
\begin{enumerate}
\item 
$$ \left\langle 
\begin{array}{c|c} 
x_1,x_2,x_3 & (x_2 x_3)^2=(x_3 x_2)^2, (x_1 x_3)^2=(x_3 x_1)^2,\\
            & [x_1, x_2]=[x_2, x_3 x_1 x_3^{-1}]=e 
\end{array}
\right\rangle$$  
\item $\langle  x_1,x_2,x_3 \  |\  (x_3 x_2 x_1)^2= (x_2 x_1 x_3)^2 = (x_1 x_3 x_2)^2 \rangle$
\item $\Z \oplus \langle x_1,x_2 | (x_1 x_2)^2=(x_2 x_1)^2 \rangle$
\item $\Z \oplus \F_2$
\item $\Z^2 \oplus \langle x_1,x_2 \ | \ (x_1 x_2)^2=(x_2 x_1)^2=e \rangle$
\end{enumerate}
\end{prop}  

The proof of Proposition \ref{th1-2} is presented completely in Section \ref{sec:1}. 
Otherwise, we skip all the braid monodromy computations, and give only the presentation 
of the fundamental group obtained by the van Kampen Theorem. These computations and the proofs 
of Propositions \ref{th1} and \ref{th2} appear in Sections \ref{sec:2} and \ref{sec:3}, 
respectively.
 
\medskip

A group is called {\it big}  if it contains a subgroup which is free (generated by 
two or more generators). By the above results, we have the following corollary:

\begin{cor}\label{bigness}
All the possibilities for fundamental groups of conic-line arrangements consist of 
two tangent conics and up to two additional lines 
are big.  
\end{cor}

The proof of this corollary appears in Section \ref{sec_big}.

\medskip

We rule out some possibilities for conic-line arrangements by the following remark.

\begin{rem}\label{remark}
The only possibility for a line to be tangent to both conics is that this line will pass through 
one of the tangency points of the two conics (fourth case in Section \ref{2conics1line_other}). 
The reason is the following. If a line was tangented to both conics, it would be a common point 
of the two duals to the two conics. Now, since the dual curve of two tangented conics is two 
tangented conics too, there are only two common points to the two dual conics. 
These two common points correspond to the two tangency points between the original two conics.  
Hence, a different line which tangents to both conics is impossible.
\end{rem}

\medskip

The paper is organized as follows. 

Section \ref{local} presents the local computations related to the 
singular points appearing in the conic-line arrangements which we deal with.

Section \ref{sec:1} deals with the proof of Proposition \ref{th1-2}.
In Section \ref{sec:2} we compute the different cases of conic-line arrangements 
which consist of two tangent conics and one additional line. 
Section \ref{sec:3} deals with  the different cases of conic-line arrangements 
which consist of two tangent conics and two additional lines. 
In Section \ref{sec_big} we show the simple proof of Corollary \ref{bigness}.

Section \ref{general} deals with some conjectures concerning the connection  
between the presentation of the fundamental group of 
the complements and the geometry of an interesting family of conic-line 
arrangements.  

In the appendix, we list all the possibilities 
for two tangent conics and two additional lines 
with the corresponding fundamental groups of their complements.

\section{Local computations}\label{local}

In this section, we present the local computations related to the 
singular points appeared  in the conic-line arrangements which we deal with.

\subsection{A singular point with two components}\label{2_comps}

In this section, we compute the local braid monodromy of a tangency point of 
two conics.

\begin{figure}[h]
\epsfysize=4.5cm  
\epsfbox{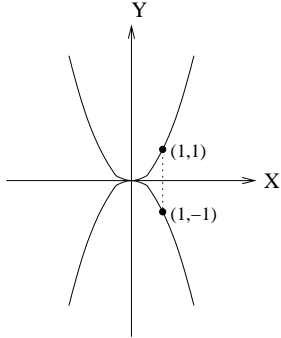}  
\caption{The singularity of $(y+x^2)(y-x^2) = 0$ at $(0,0)$}\label{2TanCon}  
\end{figure}

The tangency point between two conics can be presented locally by the 
equation: $(y+x^2)(y-x^2) = 0$ (see Figure \ref{2TanCon}).  

The singular point is $(0,0)$, and the points of  
the curve in the fiber over $x=1$ are $y=1$ and $y=-1$. 
Now, let us take a loop in $x$-axis based at $x=1$ and circumscribing $0$:   
$\alpha (t) =e^{2 \pi it}$, where  $0 \leq t \leq 1$.  Lifting it to
the curve, we get two paths on the curve:  
$\alpha _1(t) = (e^{2 \pi it}, e^{4 \pi it})$ and   
$\alpha _2(t) = (e^{2 \pi it}, -e^{4 \pi it})$ where  $0 \leq t \leq 1$.  
Projecting the paths to the $y$-axis, we get $y_1(t) = e^{4 \pi it}$ and  
 $y_2(t)= -e^{4 \pi it}$ where $0 \leq t \leq 1$.  
By substituting $t=1$, the local braid monodromy of this point 
is two counterclockwise full-twists. Hence, we have:

\begin{thm}\label{thm21}
The local braid monodromy of the tangency point of two conics is two counterclockwise
full-twists (similar to the tangency point of a conic and a line). Hence, we have 
the same induced relations: 
$$(x_1 x_2)^2 = (x_2 x_1)^2$$
where $\{ x_1,x_2 \}$ are the generators of the standard g-base.
\end{thm}

Now, we compute the corresponding Lefschetz diffeomorphism: 
\begin{cor}
The Lefschetz diffeomorphism of a tangency point between two tangent conics 
is a counterclockwise full-twist (similar to the usual tangency point).
\end{cor}

The reason for this corollary is since the Lefschetz diffeomorphism is obtained
by going on the loop mentioned above from $t=\frac{1}{2}$ 
to $t=1$. Along this interval, the two points in the fiber make one counterclockwise 
full-twist, and hence the Lefschetz diffeomorphism of this singular point
is a full-twist.

\subsection{Singular points with three components}\label{3_comps}

In this section, we  compute the local braid monodromy of four types 
of singular points consist of three components (where at most two of them 
are conics which are tangent). Then, we  compute 
the relations induced from these singular points by the 
classical van Kampen Theorem (see \cite{vK}). 

\subsubsection{First type}\label{3comps_type1}
 
\begin{figure}[h]
\epsfysize=2.5cm  
\epsfbox{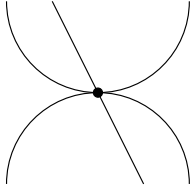}  
\caption{The singularity of $(2x+y)(y+x^2)(y-x^2) = 0$ at $(0,0)$}\label{type_b}  
\end{figure}

The local equation of the singularity of the first type is 
$(2x+y)(y+x^2)(y-x^2) = 0$ (see Figure \ref{type_b}). 

For computing the braid monodromy, we take a loop around $x=0$ and look at what happens
to the points of the curve $(2x+y)(y+x^2)(y-x^2) = 0$ in the fibers above this loop. 
 
Let $x=e^{2\pi i t}$ where $0 \leq t \leq 1$. For $t=0$ we have that $x=1$ and 
the points of the curve in the fiber over $x=1$ are $y=1,-1$ and $y=-2$. 
For $t=\frac{1}{2}$, we have $x=-1$, and the points of the curve 
in the fiber over $x=-1$ are $y=-1,1$ and $y=2$.  The point $y=1$
in the fiber $x=-1$ corresponds to the point $y=1$ in the fiber $x=1$, and similarly
the point $y=-1$ in the fiber $x=-1$ corresponds to the point $y=-1$ 
in the fiber $x=1$. The point $y=2$ in the fiber $x=-1$ corresponds to the point $y=-2$ 
in the fiber $x=1$. Hence, we get that from $t=0$ to $t=\frac{1}{2}$, the points $y=1$
and $y=-1$ do a counterclockwise full-twist and the point $y=-2$ does 
a counterclockwise half-twist around the points $y=1,-1$. 
If we continue to $t=1$, the points $y=1$
and $y=-1$ do two counterclockwise full-twists and the point $y=-2$ does 
a counterclockwise full-twist around the points $y=1,-1$ together. 

\begin{figure}[h]
\epsfysize=8cm  
\epsfbox{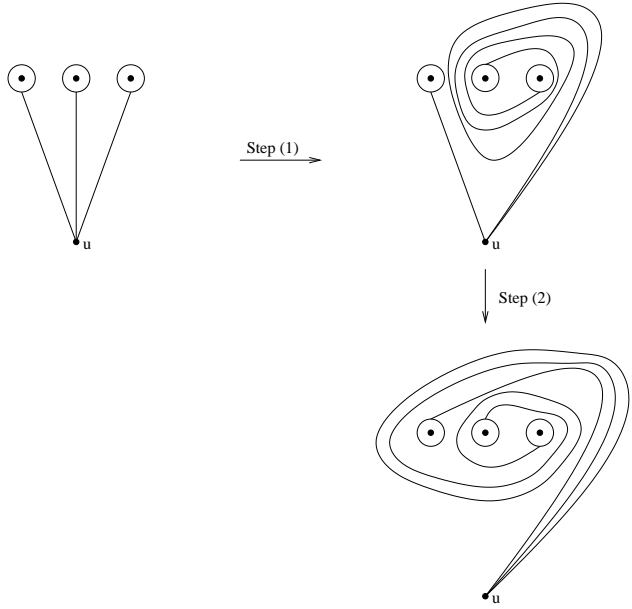}  
\caption{The g-base obtained from the standard g-base by the action of the local braid monodromy}\label{gbase_type_b}  
\end{figure}

Now, we want to compute the induced relations from the braid monodromy 
of the singular point. According to van 
Kampen's Theorem, one should compute the g-base obtained by applying the action induced
by the local braid monodromy of the singular point on the standard g-base. 
This is shown in Figure \ref{gbase_type_b}: in the first step, we perform two 
counterclockwise full-twists of points $2$ and $3$ which represent 
$y=-1$ and $y=1$ respectively. In the second step, 
we perform a counterclockwise full-twist of point $1$ (which represents 
the point $y=-2$) around the points $2$ and $3$.

Now, by van Kampen's Theorem, we get the following induced relations from 
the new g-base (where $\{ x_1,x_2,x_3 \}$ are the generators of the standard g-base):
\begin{enumerate}
\item $x_1 = x_3 x_2 x_1 x_2^{-1} x_3^{-1}$  
\item $x_2 = x_3 x_2 x_1 x_3 x_2 x_3^{-1} x_1^{-1} x_2^{-1} x_3^{-1}$
\item $x_3 = x_3 x_2 x_1 x_3 x_2 x_3 x_2^{-1} x_3^{-1} x_1^{-1} x_2^{-1} x_3^{-1}$
\end{enumerate}

From Relation (1), we get that $x_1 x_3 x_2 = x_3 x_2 x_1$. Relation (2) becomes
$x_3 x_2 x_1 x_3 x_2 = x_2 x_3 x_2 x_1 x_3$. 

Now, we show that Relation (3) is not needed. Simplifying Relation (3) yields
$x_2 x_1 x_3 x_2 x_3 = x_3 x_2 x_1 x_3 x_2$. By Relation (1), 
this relation is equal to Relation (2), and hence Relation (3) is redundant.   

To summarize, we proved the following:

\begin{thm}\label{thm23}
The local braid monodromy of the singularity presented locally by the equation 
$(2x+y)(y+x^2)(y-x^2)=0$ is: two points (corresponding to $y=1$ and $y=-1$) do two 
counterclockwise full-twists and the third point (corresponds to $y=-2$) does  
a counterclockwise full-twist around them. 

The induced relations of this point are: 
$$x_1 x_3 x_2 = x_3 x_2 x_1 \qquad ; \qquad x_3 x_2 x_1 x_3 x_2 = x_2 x_3 x_2 x_1 x_3$$
where $\{ x_1,x_2,x_3 \}$ are the generators of the standard g-base.
\end{thm}

Since the Lefschetz diffeomorphism of the singular point is obtained by computing
the action only on half of the unit circle (from $t=\frac{1}{2}$ to $t=1$), 
we have the following corollary:

\begin{cor}
The Lefschetz diffeomorphism of the singular point presented locally by
$(2x+y)(y+x^2)(y-x^2) = 0$ is a counterclockwise full-twist of the points
corresponding to $y=1$ and $y=-1$, and a half-twist of the point corresponding to $y=2$ 
(which becomes $y=-2$) around them.
\end{cor}

\begin{rem}
One can easily see that the singular point presented locally by the equation 
$y(2x+y)(y+x^2)=0$ (see Figure \ref{type_g2}) has the same braid monodromy 
and the same induced 
relations as the singular point which we have dealt with in this section, 
since these two singular points are locally the same (from the topological point 
of view).

\begin{figure}[h]
\epsfysize=2.5cm  
\epsfbox{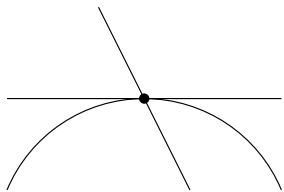}  
\caption{The singularity of $y(2x+y)(y+x^2) = 0$ at $(0,0)$}\label{type_g2}  
\end{figure}  
  
\end{rem}

\subsubsection{Second type}\label{3comps_type2}

\begin{figure}[h]
\epsfysize=2.5cm  
\epsfbox{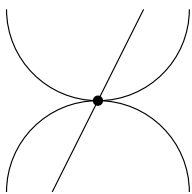}  
\caption{The singularity of $(2x-y)(y+x^2)(y-x^2) = 0$ at $(0,0)$}\label{type_c}  
\end{figure}

The local equation of the singularity of the second type is 
$(2x-y)(y+x^2)(y-x^2) = 0$  (see Figure \ref{type_c}). 

Applying almost the same computations we have performed for the first type, 
we get the following result:

\begin{thm}
The local braid monodromy of the singularity presented locally by the equation 
$(2x-y)(y+x^2)(y-x^2)=0$ is: two points (corresponding to $y=1$ and $y=-1$) do two 
counterclockwise full-twists and the third point (corresponding to $y=2$) do 
a counterclockwise full-twist around them. 

The corresponding induced relations are:    
$$x_3 x_2 x_1 = x_2 x_1 x_3 \qquad ; \qquad x_3 x_2 x_1 x_2 x_1 = x_1 x_3 x_2 x_1 x_2$$
where $\{ x_1,x_2,x_3 \}$ are the generators of the standard g-base.

The Lefschetz diffeomorphism of the singular point 
is a counterclockwise full-twist of the points
corresponding to $y=1$ and $y=-1$, and a half-twist of the point corresponding to $y=-2$ 
(which becomes $y=2$) around them.
\end{thm}

\begin{rem}

As in the previous section, one can easily see that the singular point presented 
locally by the equation $y(2x-y)(y+x^2)=0$ (see Figure \ref{type_g}) 
has the same braid monodromy and the same induced 
relations as the singular point which we have dealt with in this section, 
since these two singular points are locally the same.

\begin{figure}[h]
\epsfysize=2.5cm  
\epsfbox{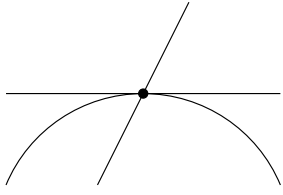}  
\caption{The singularity of $y(2x-y)(y+x^2) = 0$ at $(0,0)$}\label{type_g}  
\end{figure}  
  
\end{rem}

\subsubsection{Third type}\label{3comps_type3}

\begin{figure}[h]
\epsfysize=2.5cm  
\epsfbox{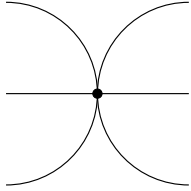}  
\caption{The singularity of $y(y^2+x)(y^2-x) = 0$ at $(0,0)$}\label{type_h}  
\end{figure}

The local equation of the singularity of the third type is $y(y^2+x)(y^2-x) = 0$ 
(see Figure \ref{type_h}). One should notice that there is a major difference between
this type of singularity and the two previous types: In this singularity there are
two ``hidden'' branch points. That is, at any fiber one has three real points 
and two complex points (i.e. complex level $2$), and in each side of the 
singularity, the complex points belong to a different conic (by the 
singularity, two real points become complex and two complex points become real).   

For computing the braid monodromy, we take a loop around $x=0$ and look what happens
to the points of the curve $y(y^2+x)(y^2-x) = 0$ in the fibers over this loop. 
 
Let $x=e^{2\pi i t}$ where $0 \leq t \leq 1$. For $t=0$ we have that $x=1$ and 
the points of the curve in the fiber over $x=1$ are $y=-1,0,1,i$ and $y=-i$. 
For $t=\frac{1}{2}$, we have $x=-1$, and we get again that the points of the curve 
in the fiber over $x=-1$ are $y=-1,0,1,i$ and $y=-i$.  By a careful checking, 
one can see that the point $y=1$
in the fiber $x=-1$ corresponds to the point $-i$ in the fiber $x=1$, and similarly
all the points except for $y=0$ in the fiber $x=-1$ 
made a $90^{\circ}$ rotation counterclockwise around $y=0 \in \C$ 
from their corresponding points in the fiber $x=1$. The point $y=0$ remains fixed. 
Now, when we continue to $t=1$, the points of the fiber continue 
to move counterclockwise around $y=0$. When we reach $t=1$ and 
reach back the fiber $x=1$, the points in the fiber made a $180^{\circ}$ rotation 
counterclockwise around $y=0 \in \C$ from their corresponding points in 
the initial fiber $x=1$. Hence, we have: 

\begin{thm}
The action of local braid monodromy 
of the point presented locally by the equation: $y(y^2+x)(y^2-x)=0$ is a 
$180^{\circ}$ rotation counterclockwise of the four points 
(corresponding to $y=1,-1,i,-i$) around the point corresponding to $y=0$, 
as shown schematically in Figure \ref{bm_type_h}.
 
\begin{figure}[h]
\epsfysize=2.5cm  
\epsfbox{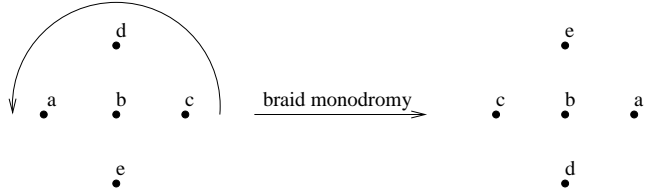}  
\caption{The local braid monodromy of the singularity $y(y^2+x)(y^2-x) = 0$}\label{bm_type_h}  
\end{figure}  
\end{thm}
 
\begin{figure}[h]
\epsfysize=12cm  
\epsfbox{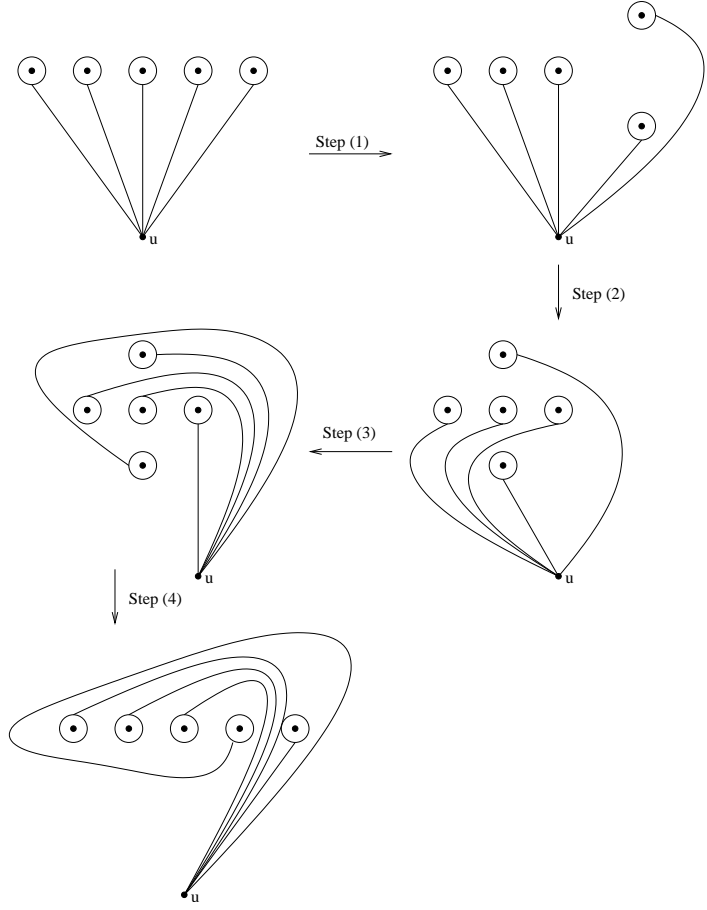}  
\caption{The g-base obtained from the standard g-base by the action of the local braid monodromy}\label{gbase_type_h}  
\end{figure}

Now, we want to compute the induced relations of this singular point. By van 
Kampen's Theorem, one should compute the g-base obtained by applying the action induced
by the local braid monodromy of the singular point on the standard g-base. Since we have 
two complex points in the fiber  before the action of the braid monodromy and after it,
we have to start by rotating the two rightmost points by $90^{\circ}$ 
counterclockwise, for representing the two complex points (see Step (1) in 
Figure \ref{gbase_type_h}). Then, we move the two complex points to be $i$ and 
$-i$ (see Step (2) there). 
Now, we apply the action of the local braid monodromy (Step (3)),
and then we return the two new complex points to the right side, and return them to the 
real axis by rotating them by $90^{\circ}$ clockwise (Step (4)). 

\begin{rem}\label{rem_clockwise}
One should notice here that when we return the points back to the real axis in Step (4), 
we rotate the points clockwise. The reason 
is that we want to compute the pure action of the braid on the g-base, and 
the clockwise rotation cancels the counterclockwise rotation which we have performed 
in Step (1).
In the global computations, we indeed rotate always counterclockwise, since we work 
there in a model (for simplifying the computations), and we can define the model 
as we want. We  get the same results in the global computations even if we rotate 
clockwise the complex points into real points.    
\end{rem}

Now, by van Kampen's Theorem, we get the following induced relations from 
the new g-base (where $\{ x_1,x_2,x_3,x_4,x_5 \}$ are the generators of the standard 
g-base):
\begin{enumerate}
\item $x_1 = x_4 x_3 x_4^{-1}$ 
\item $x_2 = x_4 x_3 x_2 x_3^{-1} x_4^{-1}$
\item $x_3 = x_4 x_3 x_2 x_1 x_2^{-1} x_3^{-1} x_4^{-1}$
\item $x_4 = x_5$
\item $x_5 = x_5 x_4 x_3 x_2 x_1 x_4 x_1^{-1} x_2^{-1} x_3^{-1} x_4^{-1} x_5^{-1}$  
\end{enumerate}

Relation (2) is equivalent to $x_4 x_3 x_2 = x_2 x_4 x_3$.
Substituting $x_1$ by $x_4 x_3 x_4^{-1}$ in Relation (3) yields $x_3 = x_4 x_3 x_2 x_4 x_3 x_4^{-1} x_2^{-1} x_3^{-1} x_4^{-1}$, which is equivalent to $x_4 x_3 x_2 x_4 x_3 = x_3 x_4 x_3 x_2 x_4$. By Relation (2), it can also be written as $x_4 x_3 x_2 x_4 x_3 = x_3 x_2 x_4 x_3 x_4$.

We  show that Relation (5) is redundant. First, we can cancel $x_5$, so we get:
$x_4 x_3 x_2 x_1 x_4 = x_5 x_4 x_3 x_2 x_1$. Since $x_4=x_5$, 
we can cancel another $x_5$,
and by substituting $x_1$ by $x_4 x_3 x_4^{-1}$ and some simplifications we get 
$x_3 x_2 x_4 x_3 x_4 = x_4 x_3 x_2 x_4 x_3$, which is equal to Relation (3). 
Hence, Relation (5) is redundant.

Therefore, we have the following corollary:

\begin{cor}
The induced relations for the singular point presented locally by the equation 
$y(y^2+x)(y^2-x)=0$ are:
\begin{enumerate}
\item $x_4 x_3 x_2 = x_2 x_4 x_3$
\item $x_3 x_2 x_4 x_3 x_4 = x_4 x_3 x_2 x_4 x_3$
\item $x_1 = x_4 x_3 x_4^{-1}$
\item $x_4 = x_5$
\end{enumerate}
where $\{ x_1,x_2,x_3,x_4,x_5 \}$ are the generators of the standard g-base.
\end{cor}

Now, we compute the corresponding Lefschetz diffeomorphism: 
\begin{cor}
The Lefschetz diffeomorphism of the singular point presented locally by 
$y(y^2+x)(y^2-x)=0$ is a $90^{\circ}$ rotation counterclockwise of the
four points around the fixed point $y=0 \in \C$. 
\end{cor}

\subsubsection{Fourth type}\label{3comps_type4}
 
This type is slightly different from the previous three types. In the previous types,
the line intersects the tangency point, but was not tangent to it. 
In this type, the additional line and the tangency point have a common tangent.

\begin{figure}[h]
\epsfysize=2.5cm  
\epsfbox{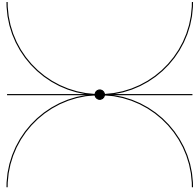}  
\caption{The singularity of $y(y+x^2)(y-x^2) = 0$ at $(0,0)$}\label{type_d}  
\end{figure}

The local equation of the singularity of this type is 
$y(y+x^2)(y-x^2) = 0$  (see Figure \ref{type_d}). 

By similar arguments as in Theorem \ref{thm21}, we get

\begin{thm}
The local braid monodromy of the singularity presented locally by the equation 
$y(y+x^2)(y-x^2)=0$ is: the three points do two 
counterclockwise generalized full-twists.

The induced relations of this point are: 
$$(x_3 x_2 x_1)^2 = (x_1 x_3 x_2)^2 = (x_2 x_1 x_3)^2,$$
where $\{ x_1,x_2,x_3 \}$ are the generators of the standard g-base.

The Lefschetz diffeomorphism is a full-twist of all the three points. 
\end{thm}

\bigskip

This result can be easily generalized to an arbitrary number 
of tangented components at the same point with a common tangent:

\begin{cor}
Let $P$ be a singular point which consists of $n$ smooth branches 
with a common tangent at $P$. Let $x_i$ be the generator which corresponds
to the $i$th branch. Then, the local braid monodromy is two generalized full-twists 
of the whole segment $[1,n]$. Moreover, the relations induced by 
this singular point are:
$$(x_n x_{n-1} \cdots x_1)^2 = (x_{n-1} \cdots x_1 x_n)^2 = \cdots = (x_1 x_n \cdots x_2)^2.$$
\end{cor}

Notice that if the tangented components have a higher order of tangency (i.e. 
the higher derivatives of the components are also equal up to the $k$th derivative,
for a given $k$), then the braid monodromy will consist of a higher power of 
a full-twist, and a higher exponent will appear in the relations.   

\subsection{Singular points with four components}\label{4_comps}

In this section, we  compute the local braid monodromy of five cases 
of singular points consist of four components, where at most two 
of them are conics. 
Then, we  compute the relations induced from these singular 
points by the classical van Kampen Theorem (see \cite{vK}). 

A singular point consists of four components in a conic-line arrangements with 
two tangent conics and up to two additional lines can be of the following 
two types: take the tangency point of the two conics, and add two lines in 
the following two ways: 
one way is to add one line which will be tangent to both conics, 
and the second line will intersect both conics (and the line) at the 
singular point. The second way will be to add two intersecting lines 
at the tangency point of the two conics.

\subsubsection{Two tangent conics with a tangent line and an intersecting line}

As before, there are three possibilities for the intersecting line: it can be
locally presented as the line $y=2x$, $y=-2x$ or $x=0$. In spite of the
fact that one can move from the first type to the second type by rotating 
the line, so locally the singularities are equivalent, 
from the global point of view 
these singularities induce different relations in the global fundamental group.
     
\medskip

\paragraph{\bf{First type}}\label{4comps_type1}

\begin{figure}[h]
\epsfysize=2.5cm  
\epsfbox{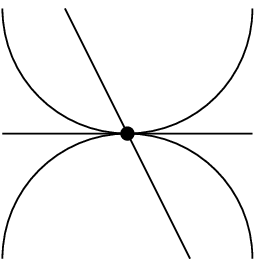}  
\caption{The singularity of $y(2x+y)(y+x^2)(y-x^2) = 0$ at $(0,0)$}\label{type_e}  
\end{figure}

The local equation of the singularity of the first type is 
$y(2x+y)(y+x^2)(y-x^2) = 0$ (see Figure \ref{type_e}). 

By a similar argument as in Theorem \ref{thm23}, we have 

\forget
For computing the braid monodromy, we take a loop around $x=0$ and 
look what happens to the points of the curve $y(2x+y)(y+x^2)(y-x^2) = 0$ 
in the fibers above this loop. 
 
Let $x=e^{2\pi i t}$ where $0 \leq t \leq 1$. For $t=0$ we have that 
$x=1$ and the points of the curve in the fiber over $x=1$ are $y=1,0,-1$ 
and $y=-2$. For $t=\frac{1}{2}$, we have $x=-1$, 
and the points of the curve in the fiber over $x=-1$ are $y=1,0,-1$ and $y=2$.  
The point $y=1$
in the fiber $x=-1$ corresponds to the point $y=1$ in the fiber $x=1$, and similarly
the point $y=-1$ in the fiber $x=-1$ corresponds to the point $y=-1$ 
in the fiber $x=1$. The point $y=0$ is fixed, and the point $y=2$ 
in the fiber $x=-1$ corresponds to the point $y=-2$ in the fiber $x=1$. 
Hence, we get that from $t=0$ to $t=\frac{1}{2}$, the points $y=-1,0$
and $y=1$ do a counterclockwise generalized full-twist 
and the point $y=-2$ does 
a counterclockwise half-twist around these three points. 
If we continue to $t=1$, the three points $y=-1,0,1$ perform 
two counterclockwise generalized full-twists and the point $y=-2$ does 
a counterclockwise full-twist around the other three points. 

For computing the induced relations, we use the van Kampen Theorem. 
Here, one should compute the g-base obtained by applying the action induced
by the local braid monodromy of the singular point on the standard g-base. 
This is shown in Figure \ref{gbase_type_e}: in the first step, we perform two 
counterclockwise generalized full-twists of points $2,3$ and $4$ which represent 
$y=-1,0$ and $y=1$ respectively. In the second step, 
we perform a counterclockwise full-twist of point $1$ (which represents 
the point $y=-2$) around the other three points.

\begin{figure}[h]
\epsfysize=8cm  
\epsfbox{gbase_type_e.eps}  
\caption{The g-base obtained from the standard g-base by the action of the local braid monodromy}\label{gbase_type_e}  
\end{figure}  

By van Kampen's Theorem, we get the following induced relations
(where $\{ x_1,x_2,x_3,x_4 \}$ are the generators of the standard g-base):
\begin{enumerate}
\item $x_1 = x_4 x_3 x_2 x_1 x_2^{-1} x_3^{-1} x_4^{-1}$  
\item $x_2 = x_4 x_3 x_2 x_1 x_4 x_3 x_2 x_3^{-1} x_4^{-1} x_1^{-1} x_2^{-1} x_3^{-1} x_4^{-1}$
\item $x_3 = x_4 x_3 x_2 x_1 x_4 x_3 x_2 x_3 x_2^{-1} x_3^{-1} x_4^{-1} x_1^{-1} x_2^{-1} x_3^{-1} x_4^{-1}$
\item $x_4 = x_4 x_3 x_2 x_1 x_4 x_3 x_2 x_4 x_2^{-1} x_3^{-1} x_4^{-1} x_1^{-1} x_2^{-1} x_3^{-1} x_4^{-1}$
\end{enumerate}

From Relation (1), we get that $x_1 x_4 x_3 x_2 = x_4 x_3 x_2 x_1$. 
Relations (2) and (3) become 
$x_4 x_3 x_2 x_1 x_4 x_3 x_2 = x_2 x_4 x_3 x_2 x_1 x_4 x_3$ and 
$x_4 x_3 x_2 x_1 x_4 x_3 x_2 x_3 = x_3 x_4 x_3 x_2 x_1 x_4 x_3 x_2$ 
respectively. By Relation (2), Relation (3) can be written as:
$$x_4 x_3 x_2 x_1 x_4 x_3 x_2 x_3 = x_3 x_2 x_4 x_3 x_2 x_1 x_4 x_3,$$ and hence
we get from Relations (2) and (3):
$$x_4 x_3 x_2 x_1 x_4 x_3 x_2 = x_3 x_2 x_4 x_3 x_2 x_1 x_4 = x_2 x_4 x_3 x_2 x_1 x_4 x_3.$$ 

We  show now  that Relation (4) is redundant. 
From Relation (4), one can get:
$$x_3 x_2 x_1 x_4 x_3 x_2 x_4 = x_4 x_3 x_2 x_1 x_4 x_3 x_2.$$
Using Relation (1), this relation is equal to Relation (3), and hence 
it is redundant. 

Hence, we have:
\forgotten

\begin{thm}
The local braid monodromy of the singularity presented locally by the equation 
$y(2x+y)(y+x^2)(y-x^2)=0$ is: three points (correspond to $y=1,0$ and $y=-1$) do two 
counterclockwise generalized full-twists and the fourth point (corresponds to $y=-2$) do 
a counterclockwise full-twist around them. 

The induced relations from this singular point are:
$$x_1 x_4 x_3 x_2 = x_4 x_3 x_2 x_1 \quad ; \quad (x_4 x_3 x_2)^2 x_1 = x_3 x_2 x_4 x_3 x_2 x_1 x_4 = x_2 x_4 x_3 x_2 x_1 x_4 x_3,$$ 
where $\{ x_1,x_2,x_3,x_4 \}$ are the generators of the standard g-base. 

The Lefschetz diffeomorphism of the singular point is a counterclockwise generalized full-twist 
of the points correspond to $y=1,0$ and $y=-1$, and a half-twist 
of the point corresponds to $y=2$ (which becomes $y=-2$).
\end{thm}

\begin{rem}
If we delete the generator $x_1$ which corresponds 
to the line that intersects the three tangent components, we get the following 
set of relations:
$$(x_4 x_3 x_2)^2 = (x_3 x_2 x_4)^2 = (x_2 x_4 x_3)^2,$$
as expected (see Section \ref{3comps_type4}).
Similarly, if we delete the generator $x_3$ which corresponds to the tangent line, 
we get the following set of relations:
$$x_1 x_4 x_2 = x_4 x_2 x_1 \quad ; \quad x_4 x_2 x_4 x_2 x_1 = x_2 x_4 x_2 x_1 x_4,$$ 
again as expected (see Section \ref{3comps_type1}).
\end{rem}

\begin{figure}[h]
\epsfysize=2.5cm  
\epsfbox{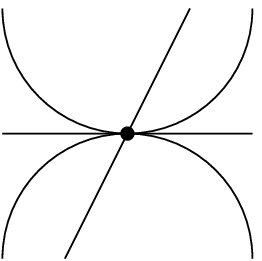}  
\caption{The singularity of $y(2x-y)(y+x^2)(y-x^2) = 0$ at $(0,0)$}\label{type_f}  
\end{figure}

\paragraph{\bf{Second type}}\label{4comps_type2}
The local equation of the singularity of the second type is 
$y(2x-y)(y+x^2)(y-x^2) = 0$  (see Figure \ref{type_f}). 

Applying almost the same computations we have performed for the previous type, 
we get:

\begin{thm}
The local braid monodromy of the singularity presented locally by the equation 
$y(2x-y)(y+x^2)(y-x^2)=0$ is: the points $y=1,0,-1$ do two 
counterclockwise generalized full-twists, and the point $y=2$ does a 
counterclockwise full-twist around the points $y=1,0,-1$ together.

The corresponding induced relations are:    
$$x_4 x_3 x_2 x_1 = x_3 x_2 x_1 x_4 \quad ; \quad (x_3 x_2 x_1)^2 x_4 = x_2 x_1 x_3 x_2 x_1 x_4 x_3 = x_1 x_3 x_2 x_1 x_4 x_3 x_2,$$ 
where $\{ x_1,x_2,x_3,x_4 \}$ are the generators of the standard g-base.

The Lefschetz diffeomorphism of the singular point presented locally by
$y(2x-y)(y+x^2)(y-x^2) = 0$ is a counterclockwise generalized full-twist 
of the points correspond to $y=1,0$ and $y=-1$, and a half-twist 
of the point corresponds to $y=-2$ (which becomes $y=2$).
\end{thm}

\begin{figure}[h]
\epsfysize=2.5cm  
\epsfbox{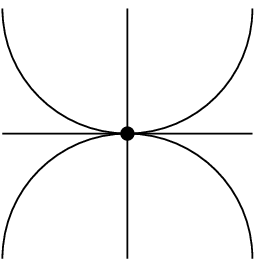}  
\caption{The singularity of $xy(y+x^2)(y-x^2) = 0$ at $(0,0)$}\label{type_j}  
\end{figure}

\paragraph{\bf{Third type}}\label{4comps_type3}
The local equation of the singularity of the third type is 
$xy(y+x^2)(y-x^2) = 0$  (see Figure \ref{type_j}). 

For computing the braid monodromy in this case, we  use the following trick:
the three components in the middle $x(y+x^2)(y-x^2)=0$ can be thought 
for a moment as a ``thick'' line which is perpendicular to the other line.
After this observation, we have two intersecting lines, 
whose local braid monodromy is a counterclockwise full-twist (even though 
we have here a thick line $x=0$, which is vertical, one can rotate it a bit 
getting two ``usual'' intersecting lines, without changing the braid monodromy).

The local braid monodromy induced by the ``thick'' line itself (which consists 
of the curve $x(y+x^2)(y-x^2)=0$) has already been computed (Section 
\ref{3comps_type3}): 
we got there that the local braid monodromy is a $180^{\circ}$ counterclockwise 
rotation of four points around one fixed point at the origin. 

Hence, the action of the local braid monodromy of the whole singular point
can be summarized as follows: 

\begin{thm}
The local braid monodromy of the singular point presented locally by the 
equation $xy(y+x^2)(y-x^2)=0$ is: first perform a $180^{\circ}$ counterclockwise 
rotation of the points correspond to $y=1,-1,i,-i$ around $y=0$, 
and then the point corresponds to $y=2$ performs a 
full-twist with a block which consists of all the other points (see Figure 
\ref{bm_type_j}).

\begin{figure}[h]
\epsfysize=2.8cm  
\epsfbox{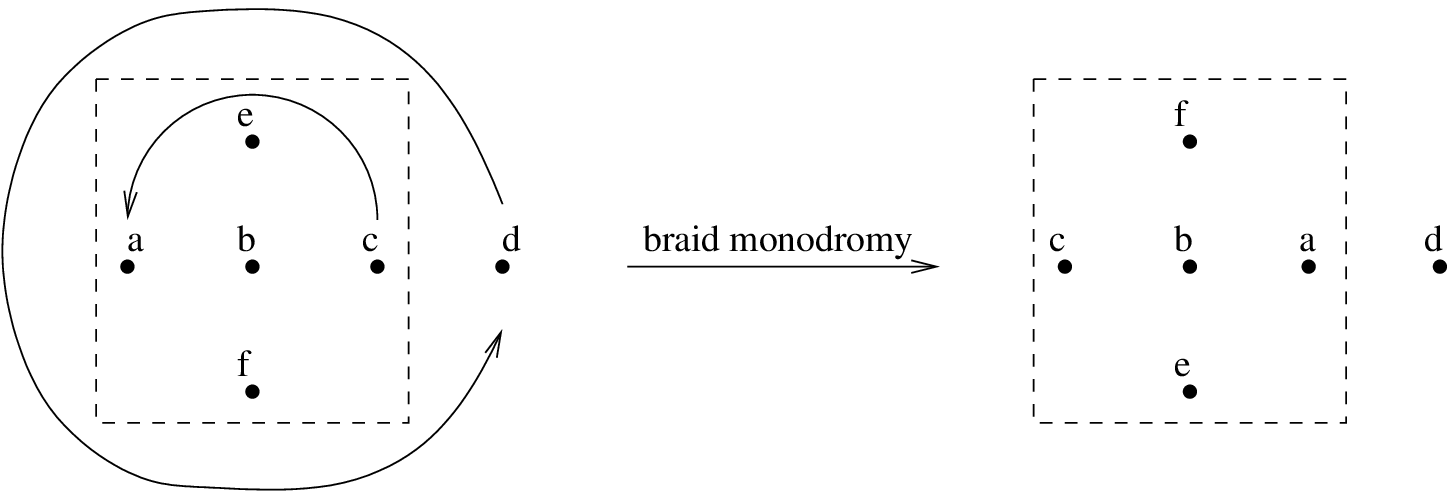}  
\caption{The action of the local braid monodromy}\label{bm_type_j}  
\end{figure}  
\end{thm}

\begin{figure}[h]
\epsfysize=12cm  
\epsfbox{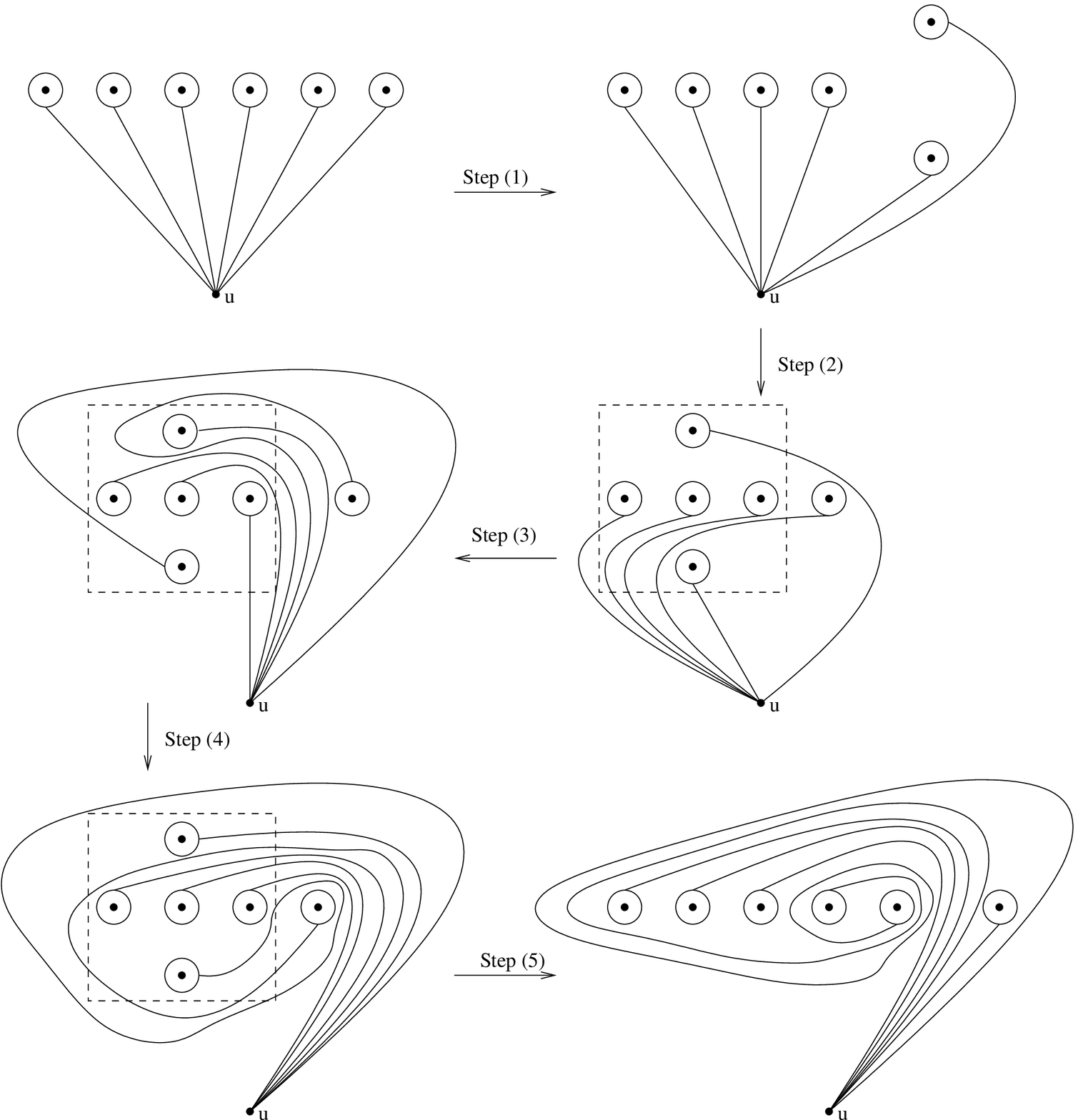}  
\caption{The g-base obtained from the standard g-base by the action of the local braid monodromy}\label{gbase_type_j}  
\end{figure}

Now, we want to compute the induced relations of this singular point. 
According to van Kampen's Theorem, one should compute the g-base 
obtained by applying the action induced by the local braid monodromy 
of the singular point on the standard g-base. Since we have 
two complex points in the fiber  before the action 
of the braid monodromy and after it, we have to start by rotating 
the two rightmost points by $90^{\circ}$ 
counterclockwise, for representing the two complex points (see Step (1) in 
Figure \ref{gbase_type_j}). Then, we move the two complex points to be over 
the second to the left real point (see Step (2) there). 
Now, we apply the action of the local braid monodromy (Steps (3) and (4)),
and then we return the two new complex points to the right side, 
and return them to the real axis by rotating them clockwise by 
$90^{\circ}$ (Step (5), see Remark \ref{rem_clockwise}). 

Now, by van Kampen's Theorem, we get the following induced relations from 
the new g-base (where $\{ x_1,x_2,x_3,x_4,x_5,x_6 \}$ are the generators of the standard g-base):
\begin{enumerate}
\item $x_1=x_5 x_4 x_3 x_4^{-1} x_5^{-1}$
\item $x_2=x_5 x_4 x_3 x_2 x_3^{-1} x_4^{-1} x_5^{-1}$
\item $x_3=x_5 x_4 x_3 x_2 x_1 x_2^{-1} x_3^{-1} x_4^{-1} x_5^{-1}$
\item $x_4=x_5 x_4 x_3 x_2 x_1 x_5 x_4 x_5^{-1} x_1^{-1} x_2^{-1} x_3^{-1} x_4^{-1} x_5^{-1}$
\item $x_5=x_6$
\item $x_6=x_6 x_5 x_4 x_3 x_2 x_1 x_5 x_4 x_5 x_4^{-1} x_5^{-1} x_1^{-1} x_2^{-1} x_3^{-1} x_4^{-1} x_5^{-1} x_6^{-1}$ 
\end{enumerate}

Relation (2) can be written as: $x_5 x_4 x_3 x_2 = x_2 x_5 x_4 x_3$.
Using Relation (1), Relation (3) becomes: 
$x_3 x_5 x_4 x_3 x_2 x_5 x_4 = x_5 x_4 x_3 x_2 x_5 x_4 x_3$.

By Relation (1) and some cancellations, Relation (4) becomes: 
$$x_4=x_5 x_4 x_3 x_2 x_5 x_4 x_3 x_4 x_3^{-1} x_4^{-1} x_5^{-1} x_2^{-1} x_3^{-1} x_4^{-1} x_5^{-1}.$$
By Relation (3), we get:
$$x_4=x_3 x_5 x_4 x_3 x_2 x_5 x_4 x_5^{-1} x_2^{-1} x_3^{-1} x_4^{-1} x_5^{-1} x_3^{-1},$$
and hence we have: $x_4 x_3 x_5 x_4 x_3 x_2 x_5 = x_3 x_5 x_4 x_3 x_2 x_5 x_4$.

Now, we  show that Relation (6) is redundant. First, this relation 
can be written as: $x_5 x_4 x_3 x_2 x_1 x_5 x_4 = x_4 x_3 x_2 x_1 x_5 x_4 x_5$. 
By Relation (1) and some cancellations, we have:
$x_5 x_4 x_3 x_2 x_5 x_4 x_3 = x_4 x_3 x_2 x_5 x_4 x_3 x_5$.
By Relation (3), we have: 
$x_3 x_5 x_4 x_3 x_2 x_5 x_4 = x_4 x_3 x_2 x_5 x_4 x_3 x_5$.
By Relation (2), this relation is equal to Relation (4), and hence Relation (6) is 
redundant.

Therefore, we get the following result for the set of relations for the singular point:

\begin{cor}
The singular point presented locally by $xy(y+x^2)(y-x^2) = 0$ 
has the following set of induced relations: 
\begin{enumerate}
\item $x_5 x_4 x_3 x_2 = x_2 x_5 x_4 x_3$
\item $x_3 x_5 x_4 x_3 x_2 x_5 x_4 = x_5 x_4 x_3 x_2 x_5 x_4 x_3 = x_4 x_3 x_5 x_4 x_3 x_2 x_5$
\item $x_1 = x_5 x_4 x_3 x_4^{-1} x_5^{-1}$
\item $x_5 = x_6$
\end{enumerate}
where $\{ x_1,x_2,x_3,x_4,x_5,x_6 \}$ are the generators of the standard g-base.
\end{cor}

Since the Lefschetz diffeomorphism of the singular point 
is obtained by computing
the action only on half of the unit circle (from $t=\frac{1}{2}$ to $t=1$), 
we have the following corollary:

\begin{cor}
The Lefschetz diffeomorphism of the singular point presented locally by
$xy(y+x^2)(y-x^2) = 0$ is a $90^{\circ}$ counterclockwise rotation of the 
points correspond to $y=1,-1,i,-i$ around $y=0$, 
and then the point corresponds to $y=2$ do a 
half-twist with a block which consists of all the other points.
\end{cor}

\subsubsection{Two tangent conics with two intersecting lines}

Although the following two types are almost a rotation of each other, we need them both, 
since the second type includes also ``hidden'' branch points inside the 
singularity, as we had in Section \ref{3comps_type3}.

\medskip

\begin{figure}[h]
\epsfysize=2.5cm  
\epsfbox{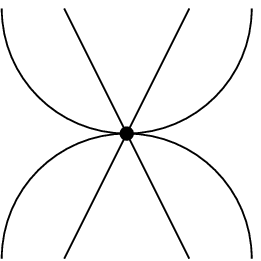}  
\caption{The singularity of $(2x+y)(2x-y)(y+x^2)(y-x^2) = 0$ at $(0,0)$}\label{type_k}  
\end{figure}

\paragraph{\bf{First type}}\label{4comps_type4}
The local equation of the singularity of the first type is 
$(2x+y)(2x-y)(y+x^2)(y-x^2) = 0$ (see Figure \ref{type_k}). 

\forget
For computing the braid monodromy, we take a loop around $x=0$ and look what happens
to the points of the curve $(2x+y)(2x-y)(y+x^2)(y-x^2) = 0$ in the fibers 
over this loop. 
 
Let $x=e^{2\pi i t}$ where $0 \leq t \leq 1$. For $t=0$ we have that $x=1$ and 
the points of the curve in the fiber over $x=1$ are $y=1,-1,2$ and $y=-2$. 
For $t=\frac{1}{2}$, we have $x=-1$, and the points of the curve 
in the fiber over $x=-1$ are $y=1,-1,-2$ and $y=2$.  The point $y=1$
in the fiber $x=-1$ corresponds to the point $y=1$ in the fiber $x=1$, and similarly
the point $y=-1$ in the fiber $x=-1$ corresponds to the point $y=-1$ 
in the fiber $x=1$. The point $y=2$ in the fiber $x=-1$ corresponds to 
the point $y=-2$ in the fiber $x=1$, and similarly the point $y=-2$ 
in the fiber $x=-1$ corresponds to the point $y=2$ in the fiber $x=1$. 
Hence, we get that from $t=0$ to $t=\frac{1}{2}$, the points $y=1$
and $y=-1$ do a counterclockwise full-twist and the points $y=2$ and $y=-2$ do 
a counterclockwise half-twist around the points $y=1,-1$. 
If we continue to $t=1$, the points $y=1$
and $y=-1$ do two counterclockwise full-twists and the points $y=2$ and $y=-2$ do 
a counterclockwise full-twist around the points $y=1,-1$ together. 

Hence we have the following result:
\forgotten

\begin{thm}
The local braid monodromy of the singularity presented locally by the equation 
$(2x+y)(2x-y)(y+x^2)(y-x^2) = 0$ is: two points (correspond to $y=1$ and $y=-1$) do two 
counterclockwise full-twists and the two other points (corresponds to $y=2$ and $y=-2$) do 
a counterclockwise full-twist around them. 

The singular point has the following induced relations:
$$x_4 x_3 x_2 x_1 = x_1 x_4 x_3 x_2 = x_3 x_2 x_1 x_4 \quad ; \quad x_4 x_3 x_2 x_1 x_3 x_2 = x_2 x_4 x_3 x_2 x_1 x_3$$  
where $\{ x_1,x_2,x_3,x_4 \}$ are the generators of the standard g-base.

The Lefschetz diffeomorphism of the singular point is a counterclockwise full-twist
of the points $y=1$ and $y=-1$, and a counterclockwise half-twist of 
the points $y=2$ and $y=-2$ around the points $y=1,-1$. 
\end{thm}

\paragraph{\bf{Second type}}\label{4comps_type5}

\begin{figure}[h]
\epsfysize=2.5cm  
\epsfbox{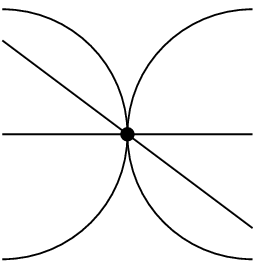}  
\caption{The singularity of $y(x+2y)(y^2+x)(y^2-x) = 0$ at $(0,0)$}\label{type_i}  
\end{figure}

The local equation of the singularity of the second type is 
$y(x+2y)(y^2+x)(y^2-x) = 0$ (see Figure \ref{type_i}). 
 
There is a major difference 
between this type of singularity and the previous one: In this singularity there are
two ``hidden'' branch points too (as in Section \ref{3comps_type3}). That is, at any fiber, one has four real points 
and two complex points (i.e. complex level $2$), and in each side of the 
singularity, the complex points belong to a different conic (since by the 
singularity, two real points become complex and two complex points become real).   

For computing the braid monodromy in this case, we  use the same trick we
have already used in Section \ref{4comps_type3}:
the two lines in the middle, $x+2y=0$ and $y=0$, can be thought for a moment as a 
``thick'' line which is perpendicular to the tangent direction of the two 
tangent conics.
We have already computed this case (Section \ref{3comps_type3}): we got there that the braid monodromy 
is a $180^{\circ}$ counterclockwise rotation of four points around one fixed point 
at the origin. 
After observing this, we should add into account that the ``thick'' line 
stands for two lines. Hence, the fixed point at the origin is now decomposed 
into two close points, which are doing a counterclockwise full-twist 
(like a usual node).  To summarize, we have: 

\begin{thm}
The local braid monodromy of the singularity presented locally 
by the equation $y(x+2y)(y^2+x)(y^2-x)=0$ is: The four points correspond to $y=1,-1,i$ and 
$y=-i$ do a $180^{\circ}$ counterclockwise rotation, and the two points in the center 
do a counterclockwise full-twist (see Figure \ref{bm_type_i}).

\begin{figure}[h]
\epsfysize=2.8cm  
\epsfbox{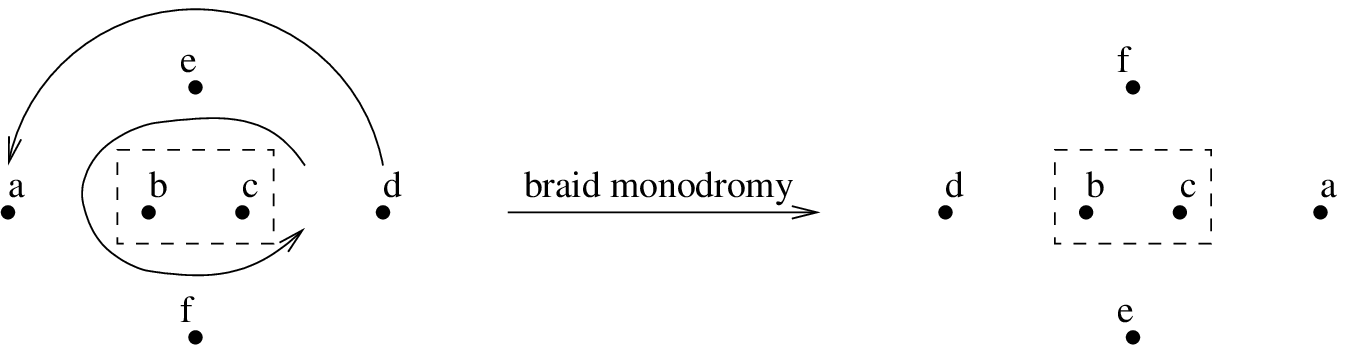}  
\caption{The action of the local braid monodromy}\label{bm_type_i}  
\end{figure}  

\end{thm}
   
\begin{figure}[h]
\epsfysize=12cm  
\epsfbox{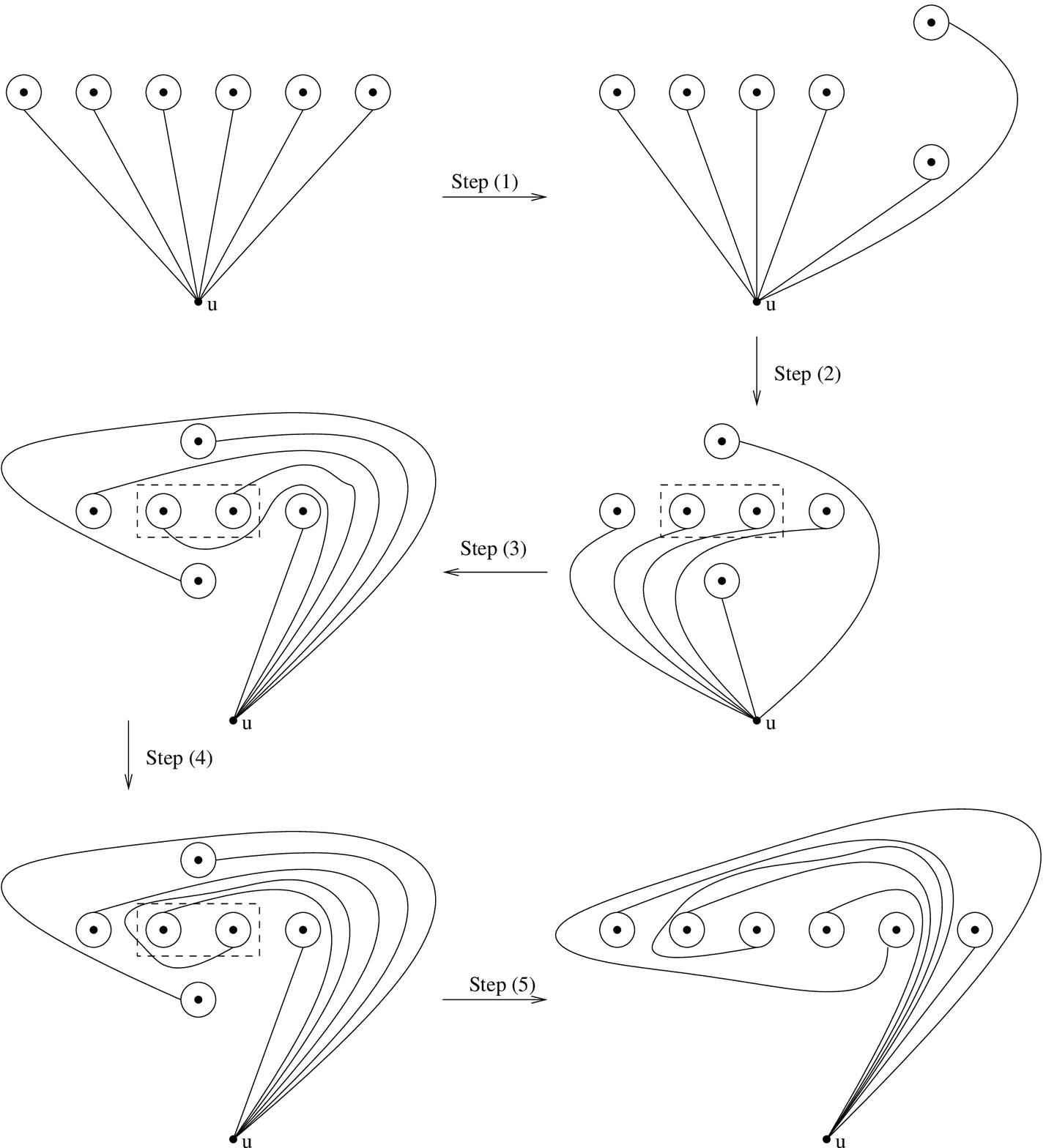}  
\caption{The g-base obtained from the standard g-base by the action of the local braid monodromy}\label{gbase_type_i}  
\end{figure}

Now, we want to compute the induced relations of this singular point. 
By van Kampen's Theorem, we compute the g-base obtained by applying the action induced
by the local braid monodromy of the singular point on the standard g-base. Since we have 
two complex points in the fiber  before the action of the braid monodromy and after it,
we have to start by rotating the two rightmost points by $90^{\circ}$ 
counterclockwise, for representing the two complex points (see Step (1) in 
Figure \ref{gbase_type_i}). Then, we move the two complex points to be over 
the two middle real points (see Step (2) there). 
Now, we apply the action of the local braid monodromy (Steps (3) and (4)),
and then we return the two new complex points to the right side, and return them to the 
real axis by rotating them by $90^{\circ}$ clockwise (Step (5), see Remark \ref{rem_clockwise}). 

Now, by van Kampen's Theorem, we get the following induced relations from 
the new g-base (where $\{ x_1,x_2,x_3,x_4,x_5,x_6 \}$ are the generators of the standard g-base):
\begin{enumerate}
\item $x_1=x_5 x_4 x_5^{-1}$
\item $x_2=x_5 x_4 x_3 x_2 x_3^{-1} x_4^{-1} x_5^{-1}$
\item $x_3=x_5 x_4 x_3 x_2 x_3 x_2^{-1} x_3^{-1} x_4^{-1} x_5^{-1}$
\item $x_4=x_5 x_4 x_3 x_2 x_1 x_2^{-1} x_3^{-1} x_4^{-1} x_5^{-1}$
\item $x_5=x_6$
\item $x_6=x_6 x_5 x_4 x_3 x_2 x_1 x_5 x_1^{-1} x_2^{-1} x_3^{-1} x_4^{-1} x_5^{-1} x_6^{-1}$ 
\end{enumerate}

From Relation (2), we get $x_5 x_4 x_3 x_2 = x_2 x_5 x_4 x_3$. Using this relation,
Relation (3) becomes: $x_3=x_2 x_5 x_4 x_3 x_4^{-1} x_5^{-1} x_2^{-1}$, and hence
$x_2 x_5 x_4 x_3 = x_3 x_2 x_5 x_4$. Using Relation (2) again, we get 
$x_5 x_4 x_3 x_2 = x_3 x_2 x_5 x_4$.

By Relation (1), Relation (4) can be written: 
$x_4 x_5 x_4 x_3 x_2 x_5 = x_5 x_4 x_3 x_2 x_5 x_4$     

Now, we  show that Relation (6) is redundant. Since $x_5=x_6$, we can simplify 
Relation (6) to the following form: $x_5 x_4 x_3 x_2 x_1 = x_4 x_3 x_2 x_1 x_5$. 
Now, by Relation (1), we get: 
$x_5 x_4 x_3 x_2 x_5 x_4 x_5^{-1} = x_4 x_3 x_2 x_5 x_4 x_5^{-1} x_5$. 
By some simplifications and Relation (3), we get Relation (4), 
and hence Relation (6) is redundant.   

Therefore, we get the following result: 
\begin{cor}
The singular point presented locally by the equation $y(x+2y)(y^2+x)(y^2-x)=0$ has 
the following set of relations:
\begin{enumerate}
\item $x_5 x_4 x_3 x_2 = x_2 x_5 x_4 x_3 = x_3 x_2 x_5 x_4$
\item $x_4 x_5 x_4 x_3 x_2 x_5 = x_5 x_4 x_3 x_2 x_5 x_4$
\item $x_1 = x_5 x_4 x_5^{-1}$
\item $x_5 = x_6$
\end{enumerate}
where $\{ x_1,x_2,x_3,x_4,x_5,x_6 \}$ are the generators of the standard g-base.
\end{cor}

If we delete the generators $x_2$ or $x_3$ which correspond 
to the lines $y=0$ or $x+2y=0$ respectively, we get the set of relations
as in Section \ref{3comps_type3} as expected.

\medskip

Since the Lefschetz diffeomorphism is a half of the action of the braid monodromy, 
we have that:  
\begin{cor}
The Lefschetz diffeomorphism of the singular point presented locally by 
$y(x+2y)(y^2+x)(y^2-x)=0$ is a counterclockwise $90^{\circ}$ rotation of 
the four points around the center (which consists of two points), and then a counterclockwise half-twist 
of the two points in the center. 
\end{cor}

\section{Two conics which are tangent to each other at two points}\label{sec:1}  

In this section, we prove Proposition \ref{th1-2} which states that if $S$ is a curve 
in $\cpt$ which composed of two tangent conics, then:  
$$\pcpt \cong \langle x_1,x_2\ |\ (x_1 x_2)^2=(x_2 x_1)^2=e \rangle$$  
  
\begin{proof}[Proof of Proposition \ref{th1-2}]  
Figure \ref{fig} shows a curve composed of two tangent conics.   
  
\begin{figure}[h]
\epsfysize=3cm  
\centerline{\epsfbox{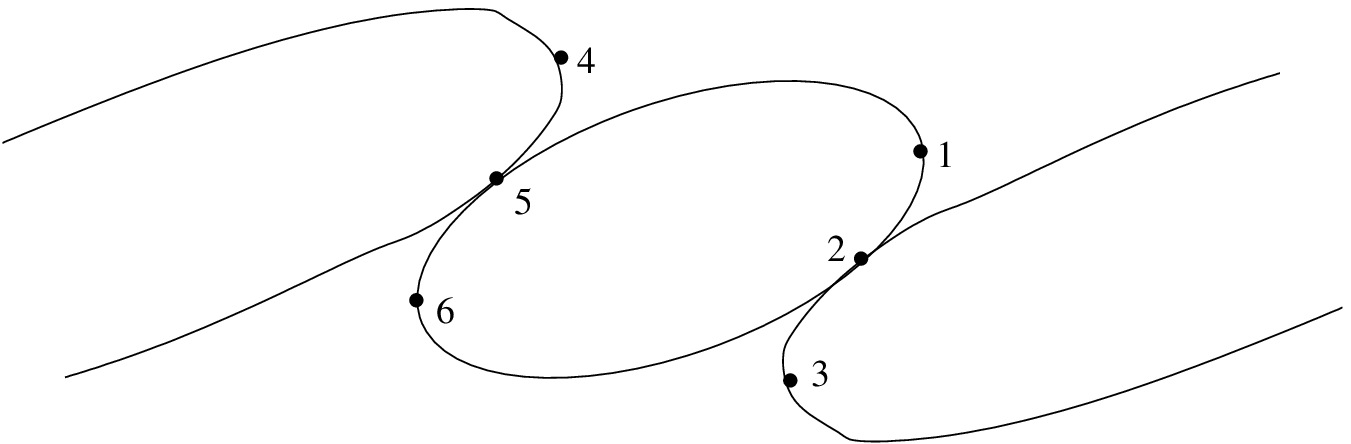}}  
\caption{Two tangent conics}\label{fig}  
\end{figure}  

For computing the braid monodromy of this curve, we first have to compute the 
Lefschetz pairs of the singular points and to identify their types for 
applying the Moishezon-Teicher algorithm correctly (see \cite{MoTe2}). 
In the following table, we summarize this data.

\begin{center}
\begin{tabular}{|c|c|c|c|}
\hline
j & $\la _{x_j}$ & $\eps _{x_j}$ & $\de _{x_j}$ \\
\hline
1 & $P_3$   & 1 & $\De^{\frac{1}{2}}_{\R I_2}<3>$ \\
2 & $<2,3>$ & 4 & $\De^2 <2,3>$ \\
3 & $<1,2>$ & 1 & $\De^{\frac{1}{2}}_{I_2 \R} <1>$ \\
4 & $P_3$   & 1 & $\De^{\frac{1}{2}} _{\R I_2}<3>$ \\
5 & $<2,3>$ & 4 & $\De^2 <2,3>$ \\
6 & $<1,2>$ & 1 & $\De^{\frac{1}{2}}_{I_2 \R} <1>$ \\
\hline
\end{tabular}
\end{center}

By the Moishezon-Teicher algorithm, we get the skeletons related to the braid monodromy
as shown in Figure \ref{thm1-bm}.

\begin{figure}[h]
\epsfysize=8cm  
\centerline{\epsfbox{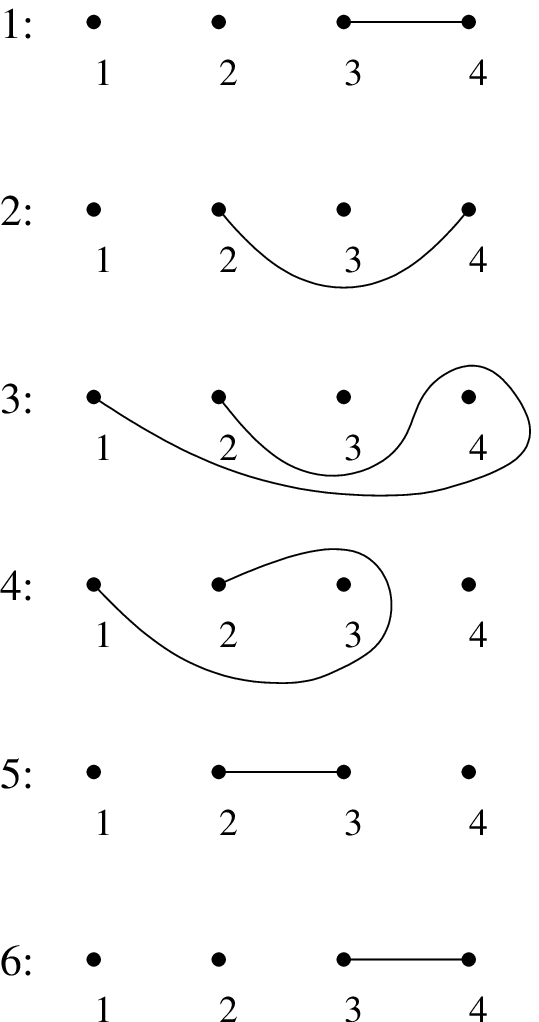}}  
\caption{Skeletons for the braid monodromy of the curve in Proposition \ref{th1-2} (Figure 
\ref{fig})}\label{thm1-bm}  
\end{figure}  
  
By the van Kampen Theorem, we get the following presentation for the group $\pi_1 (\cpt-S)$: 

\medskip

\noindent
Generators: $\{ x_1,x_2,x_3,x_4 \}$. \\  
Relations:  
\begin{enumerate}    
\item $x_4 x_3 x_2 x_1 = e$ (projective relation)  
\item $x_3 = x_4$  
\item $(x_2 x_4)^2 = (x_4 x_2)^2$  
\item $x_1 = x_4 x_2 x_4^{-1}$  
\item $x_1 = x_3 x_2 x_3^{-1}$  
\item $(x_2 x_3)^2 = (x_3 x_2)^2$  
\item $x_3 = x_4$  
\end{enumerate}    

\medskip
  
It remains to show that this presentation is equivalent to the 
presentation of the group in the formulation of the proposition.   
  
By Relation (2) (which is equal to Relation (7)), Relation (6) is equal to Relation (3)
and Relation (5) is equal to Relation (4). Hence, Relations (5),(6) and (7) are redundant.

On the other hand, the first relation becomes $x_3^2 x_2 x_1 = e$. By Relation (5), we have
$x_3^2 x_2 x_3 x_2 x_3^{-1} = e$, which yields that: $(x_3 x_2)^2=e$. We also get that 
$x_1$ is redundant, and therefore we get the requested presentation. 
\end{proof}

\section{Two tangent conics with an additional line and Proposition \ref{th1}}\label{sec:2}  

A {\it simple tangency point} can be presented locally as a tangency point between 
two smooth branches of the curve.  
We split our computations into two type of arrangements: arrangements with only 
simple tangency points and arrangements with other singular points. 

In the first subsection, we focus on arrangements with only simple tangency
points, and then we deal with the other arrangement (second subsection).

\subsection{An arrangement with a simple tangency point}\label{2conics_1line}  
  
\begin{figure}[h]
\epsfysize=4cm  
\centerline{\epsfbox{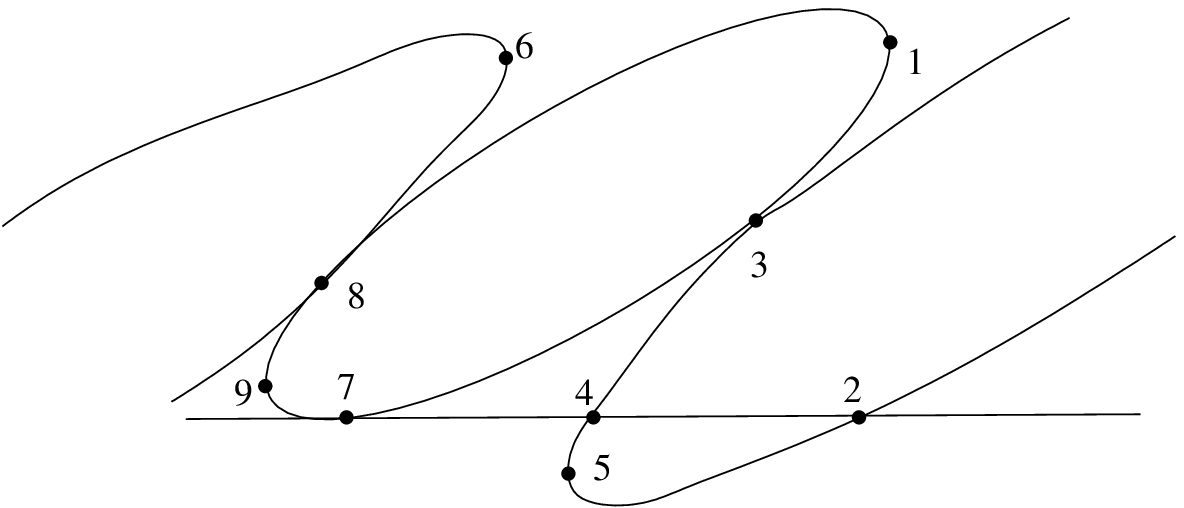}}  
\caption{The conic-line arrangement with two tangent conics and a line which tangents in a simple tangency point}\label{thm2-fig}  
\end{figure}

In this section we prove that if $S$ is a curve 
in $\cpt$ which composed of two tangent conics and one additional line which tangents to one of the conics and 
intersects the other (see Figure \ref{thm2-fig}), then:
$$\pcpt \cong \langle x_1,x_2\ |\ (x_1 x_2)^2=(x_2 x_1)^2 \rangle.$$  
This configuration is the unique configuration 
of two tangent conics with an additional tangent line, 
which tangents to the conics in 
a simple tangency point (due to Remark \ref{remark}).  

By the braid monodromy techniques and the van Kampen Theorem, we get the following 
presentation for $\pi_1(\cpt- S)$:

\medskip

\noindent  
Generators: $\{ x_1,x_2,x_3,x_4,x_5 \}$. \\  
Relations:  
\begin{enumerate}    
\item $x_5 x_4 x_3 x_2 x_1 = e$ (projective relation)  
\item $x_4= x_5$
\item $x_1 x_2 = x_2 x_1$
\item $(x_3 x_5)^2 = (x_5 x_3)^2$
\item $(x_2 x_1 x_2^{-1})(x_5 x_3 x_5^{-1}) = (x_5 x_3 x_5^{-1})(x_2 x_1 x_2^{-1})$
\item $x_2 = x_5 x_3 x_5^{-1}$
\item $x_1^{-1} x_2 x_1 = x_4 x_3 x_4^{-1}$
\item $(x_2 x_1 x_2^{-1} x_5)^2 = (x_5 x_2 x_1 x_2^{-1})^2$
\item $(x_3 x_4)^2 = (x_4 x_3)^2$
\item $x_5 = x_2 x_1^{-1} x_2^{-1} x_4 x_2 x_1 x_2^{-1}$ 
\end{enumerate}    

\medskip

By Relation (6), Relation (5) has the following form: 
$(x_2 x_1 x_2^{-1})x_2  = x_2 (x_2 x_1 x_2^{-1})$. This is equal to 
$x_1 x_2 = x_2 x_1$, which is already known by Relation (3). 
Hence, Relation (5) is redundant.

By Relation (3), Relation (7) is equal to $x_2 = x_4 x_3 x_4^{-1}$,
which is equivalent to Relation (6) by Relation (2). Hence, Relation (7) is 
redundant too. 

By Relation (3) again, Relation (8) gets the form $(x_1 x_5)^2 = (x_5 x_1 )^2$,
and Relation (10) is reduced to: $x_5 = x_1^{-1} x_4 x_1$. Since $x_4=x_5$, 
we have $x_1 x_5 =x_5 x_1$. Hence, Relation (8) is redundant.

Relation (4) and Relation (9) are equal, since $x_4=x_5$.

Therefore, we have the following equivalent presentation:

\medskip

\noindent
Generators: $\{ x_1,x_2,x_3,x_5 \}$. \\  
Relations:  
\begin{enumerate}    
\item $x_5^2 x_3 x_2 x_1 = e$
\item $x_1 x_2 = x_2 x_1$
\item $(x_3 x_5)^2 = (x_5 x_3)^2$
\item $x_2 = x_5 x_3 x_5^{-1}$
\item $x_1 x_5 = x_5 x_1$ 
\end{enumerate}    

\medskip

Substituting $x_5 x_3 x_5^{-1}$ for $x_2$ in Relation (1) and Relation (2), yields
the relations $x_5 x_3 x_5 x_3 x_1 = e$ and $x_1 x_3 = x_3 x_1$ respectively. 
Hence, we get the following presentation:

\medskip

\noindent
Generators: $\{ x_1,x_3,x_5 \}$. \\  
Relations:  
\begin{enumerate}    
\item $x_5 x_3 x_5 x_3 x_1 = e$
\item $x_1 x_3 = x_3 x_1$
\item $(x_3 x_5)^2 = (x_5 x_3)^2$
\item $x_1 x_5 = x_5 x_1$ 
\end{enumerate}    

\medskip

By the first relation, $x_1 = (x_5 x_3)^{-2}$. Substituting it for $x_1$ in Relations
(2) and (4) yields the relation $(x_5 x_3)^2 = (x_3 x_5)^2$ twice, which is 
already known by Relation (3). Hence, Relations (2) and (4) are redundant.

So the final presentation is 
$$\langle x_3,x_5 \ | \ (x_3 x_5)^2 = (x_5 x_3)^2 \rangle$$
which appears in Proposition \ref{th1} as needed. 

\begin{rem}
This fundamental group was computed independently by Degtyarev too (see
\cite[Section 3.3.5]{Deg}). 
\end{rem}

\subsection{Arrangements with other singular points}\label{2conics1line_other}

Apart from the previous case, 
we have four more possibilities to locate a line into 
a configuration of two tangent conics:
\begin{enumerate}
\item The line intersects the conic transversally. 
\item The line passes through one of the tangency points between 
the two conics, but it is not tangent to the conics 
(see Figure \ref{case1_2}).

\begin{figure}[h]
\epsfysize=4cm  
\centerline{\epsfbox{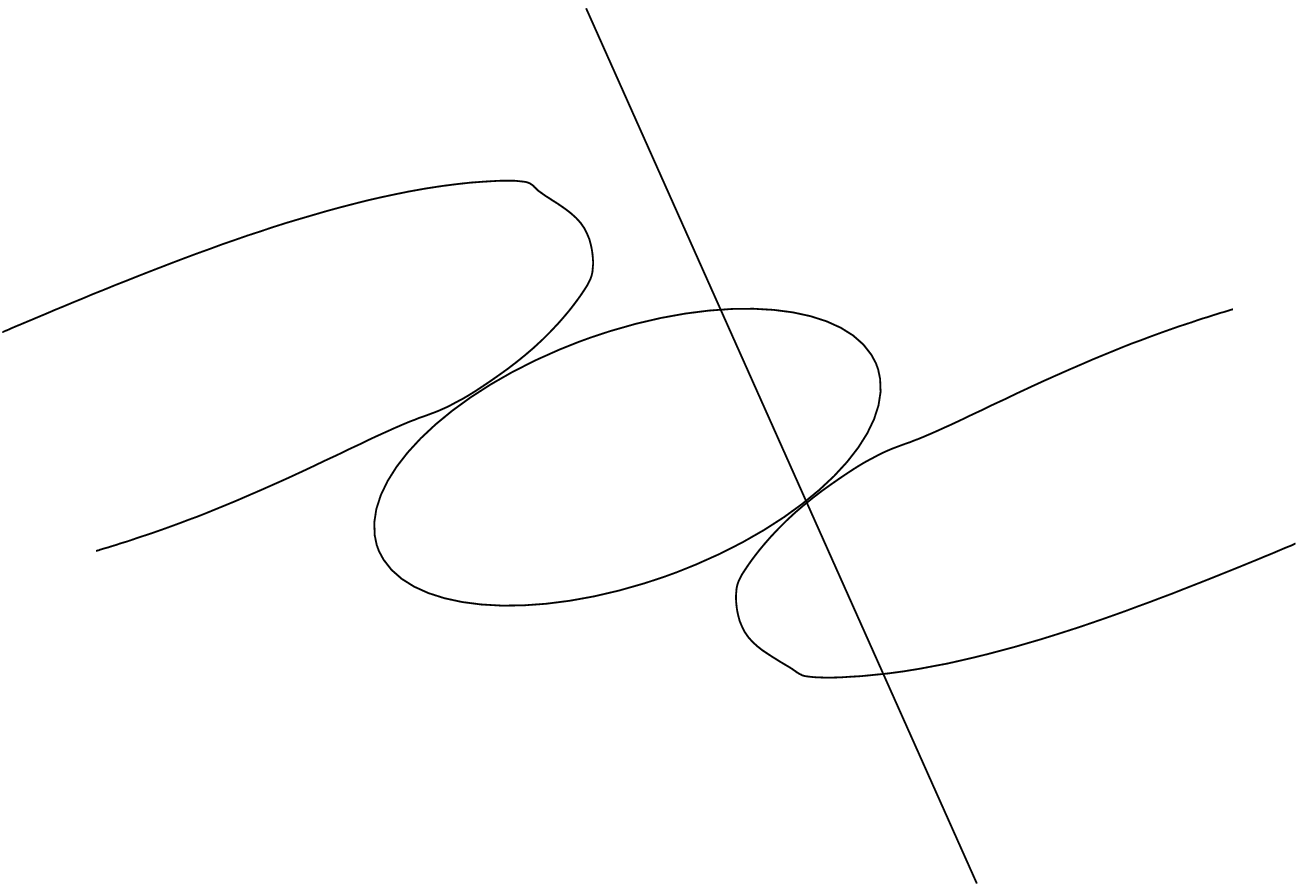}}  
\caption{The line passes through the tangency point but it is not tangent}\label{case1_2}  
\end{figure}  

\item The line passes through one of the tangency points between the two conics,
and it is tangent to the conics (see Figure \ref{case1_3}).
\begin{figure}[h]
\epsfysize=4cm  
\centerline{\epsfbox{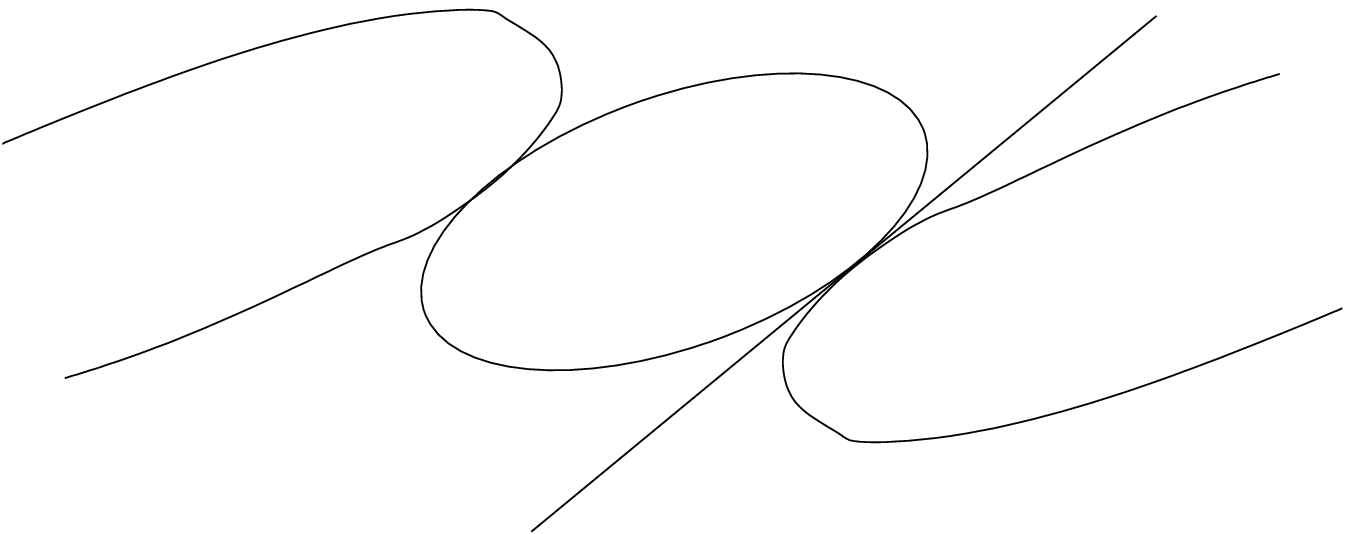}}  
\caption{The line passes through the tangency point and it is tangent}\label{case1_3}  
\end{figure}  

\item The line passes through the two tangency points (see Figure \ref{case1_1}).
\begin{figure}[h]
\epsfysize=4cm  
\centerline{\epsfbox{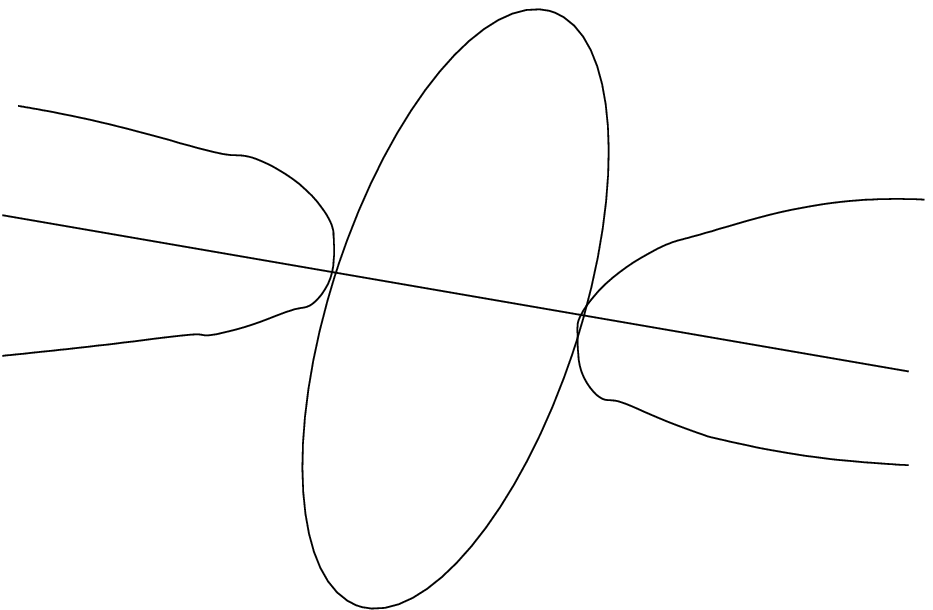}}  
\caption{The line passes through the two tangency points}\label{case1_1}  
\end{figure}  

\end{enumerate}

Based on Proposition \ref{th1-2} and \cite{Ga}, 
the fundamental group of the complement of the first case is  
$$\Z \oplus \langle x_1,x_2 \ | \ (x_1 x_2)^2=(x_2 x_1)^2=e \rangle,$$
where $e$ is the identity element of the group.

In the following subsections, we  compute the third and 
the fourth cases, 
whence the result of the second case turns out to be isomorphic to the third case.

\subsubsection{Third case}

Using braid monodromy techniques and the van Kampen theorem, 
we get the following presentation for the fundamental group 
of the complement of the curve which appear in Figure \ref{case1_3}:

\medskip

\noindent  
Generators: $\{ x_1,x_2,x_3,x_4,x_5 \}$. \\  
Relations:  
\begin{enumerate}    
\item $x_5 x_4 x_3 x_2 x_1 = e$ (projective relation)  
\item $x_4= x_5$
\item $(x_5 x_3 x_2)^2 = (x_3 x_2 x_5)^2 = (x_2 x_5 x_3)^2$
\item $x_1 = x_5 x_3 x_2 x_3^{-1} x_5^{-1}$
\item $x_1 = x_4 x_2 x_4^{-1}$
\item $(x_2 x_4)^2 = (x_4 x_2)^2$
\item $x_3^{-1} x_4 x_3 = x_5$
\end{enumerate}    

\medskip

By Relations (2) and (5), Relation (4) becomes  
$x_4 x_3 x_2 x_3^{-1} x_4^{-1} = x_4 x_2 x_4^{-1}$, and hence $x_2 x_3 =x_3 x_2$
 
By Relation (2), we get from Relation (6) that $x_3 x_4 = x_4 x_3$.

Using Relations (2) and (5) again, the first relation becomes:
$x_4 x_4 x_3 x_2 x_4 x_2 x_4^{-1}=e$. By some cancellations and Relation (6), 
we get: $x_3 = (x_2 x_4)^{-2}$. 
 
Hence, we get the following equivalent presentation:

\medskip

\noindent  
Generators: $\{ x_2,x_3,x_4 \}$. \\  
Relations:  
\begin{enumerate}    
\item $x_3 = (x_2 x_4)^{-2}$
\item $(x_4 x_3 x_2)^2 = (x_3 x_2 x_4)^2 = (x_2 x_4 x_3)^2$
\item $x_2 x_3 = x_3 x_2$
\item $(x_2 x_4)^2 = (x_4 x_2)^2$
\item $x_4 x_3 = x_3 x_4$
\end{enumerate}    

\medskip

By Relation (1) and (4), Relations (3),(5) and the right equation of Relation (2) 
become trivial.  

For the left equation of Relation (2), be Relation (1) we substitute $x_3$ by 
$(x_2 x_4)^{-2}$, to get: 
$$x_4 (x_2 x_4)^{-2} x_2 x_4 (x_2 x_4)^{-2} x_2 = (x_2 x_4)^{-2} x_2 x_4 (x_2 x_4)^{-2} x_2 x_4.$$
By cancellations, we get:
$$x_4 (x_2 x_4)^{-3} x_2 = (x_2 x_4)^{-2},$$
which becomes trivial by Relation (4) and hence it is redundant.

Therefore, we get the following final presentation:
$$\pi_1(\C\P^2-S) \cong \langle x_2,x_4 \ | \ (x_2 x_4)^2 = (x_4 x_2)^2 \rangle$$

\subsubsection{Fourth case}

Using braid monodromy techniques and the van Kampen theorem, 
we get the following presentation for the fundamental group 
of the complement of the curve which appear in Figure \ref{case1_1}:

\medskip

\noindent  
Generators: $\{ x_1,x_2,x_3,x_4,x_5 \}$. \\  
Relations:  
\begin{enumerate}    
\item $x_5 x_4 x_3 x_2 x_1 = e$ (projective relation)  
\item $x_4 = x_5$
\item $x_1 = x_4 x_3 x_4^{-1}$
\item $x_2 x_4 x_3 = x_4 x_3 x_2$
\item $x_4 x_3 x_2 x_4 x_3 = x_3 x_2 x_4 x_3 x_4$
\item $x_1 = x_2^{-1} x_3^{-1} x_4^{-1} x_3 x_4 x_3 x_2$
\item $x_2^{-1} x_3^{-1} x_4 x_3 x_2 = x_1 x_5 x_1^{-1}$
\item $x_2 x_1 x_5 = x_1 x_5 x_2$
\item $x_1 x_5 x_2 x_1 x_5 = x_5 x_2 x_1 x_5 x_1$ 
\end{enumerate}    

\medskip

By Relations (2) and (3) we replace $x_5$ and $x_1$ by $x_4$ and $x_4 x_3 x_4^{-1}$
respectively. By these replacements, Relations (4) and (8) become equal. Moreover, 
Relations (5) and (9) are the same (using Relation (4)). 
So we get the following presentation: 

\noindent  
Generators: $\{ x_1,x_2,x_3,x_4,x_5 \}$. \\  
Relations:  
\begin{enumerate}    
\item $x_4^2 x_3 x_2 x_4 x_3 x_4^{-1}= e$   
\item $x_2 x_4 x_3 = x_4 x_3 x_2$
\item $x_4 x_3 x_2 x_4 x_3 = x_3 x_2 x_4 x_3 x_4$
\item $x_4 x_3 x_4^{-1} = x_2^{-1} x_3^{-1} x_4^{-1} x_3 x_4 x_3 x_2$
\item $x_2^{-1} x_3^{-1} x_4 x_3 x_2 = x_4 x_3 x_4 x_3^{-1} x_4^{-1}$
\end{enumerate}    

Relation (1) can be simplified (using Relation (2)) to the following form: 
$x_2 x_4 x_3 x_4 x_3 = e$. Hence, we have: $x_2 =(x_4 x_3)^{-2}$.
Using this, all the other relations become trivial (by simple computations).

Hence, the resulting group is the free group with two generators.   

\begin{rem}
The above two results were computed independently by Degtyarev too (see \cite{Deg}).
\end{rem}

\section{Two tangent conics with two additional 
lines and Proposition \ref{th2}}\label{sec:3}  

As in the previous section, we divide our treatment into two cases: 
We first focus on arrangements with only simple tangency
points (first subsection), and then we deal with the other arrangement (second subsection).

\subsection{Arrangements with simple tangency points} \label{2conics_2lines} 

In this section, we present the fundamental groups of the 
complement of the two different possibilities for curves with   
two tangent conics and two additional lines 
which are tangent in only 
simple tangency points: the case where each line is tangent a different conic, 
and the case where the two lines are tangent to the same conic.

\begin{figure}[h]
\epsfysize=5cm  
\centerline{\epsfbox{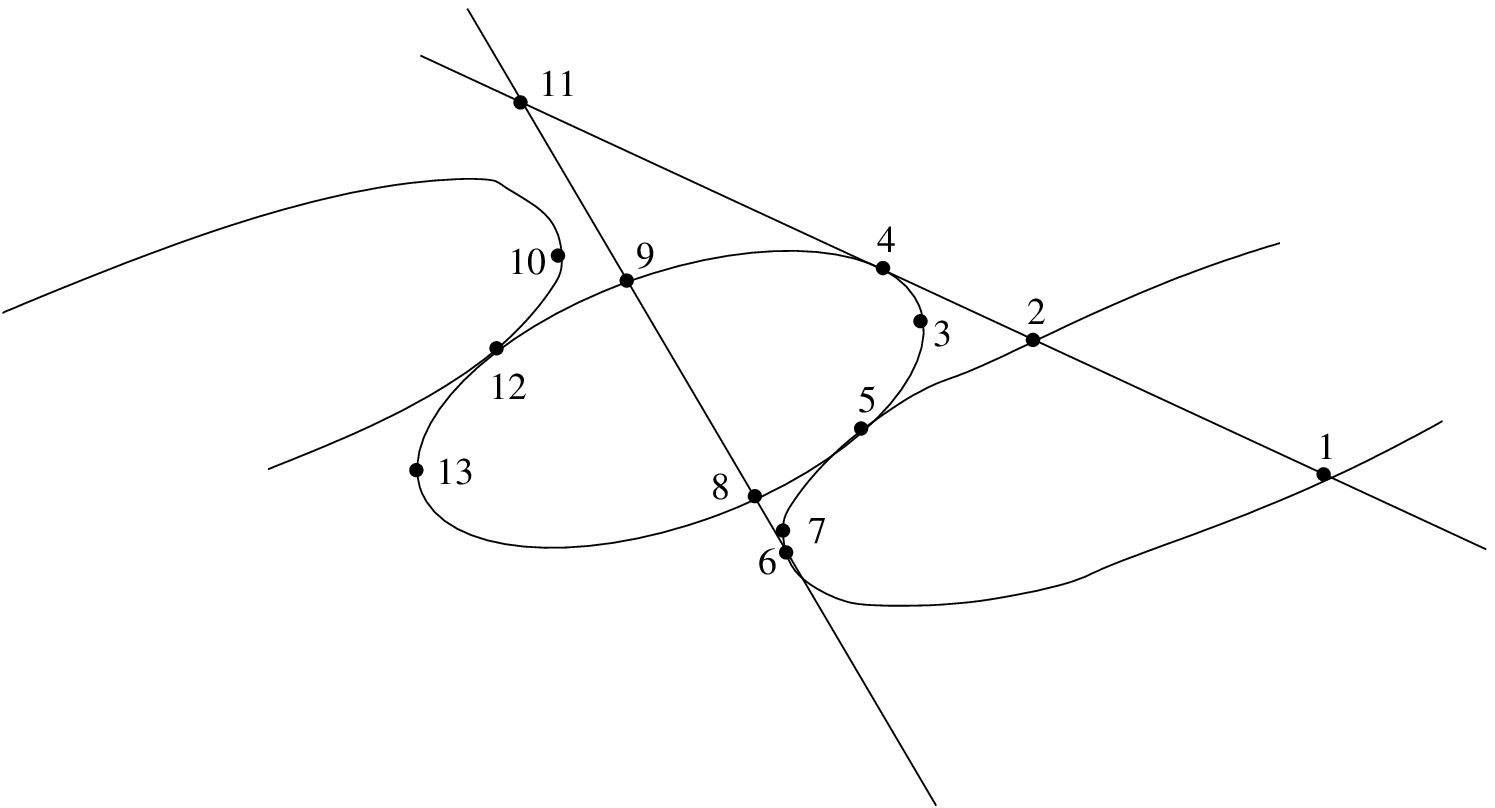}}  
\caption{The conic-line arrangement where each line is tangent to a different conic}\label{thm3-fig}  
\end{figure}

\begin{thm}
Let $S$ be a curve in $\cpt$ composed of two tangent conics 
and two additional lines, where each line is tangent to a different conic
(see Figure \ref{thm3-fig}), then:
$$\pcpt \cong \langle x_1 \rangle \oplus 
\langle x_2,x_3 \ |  \ (x_2 x_3)^2=(x_3 x_2)^2 \rangle$$
\end{thm}

\begin{figure}[h]
\epsfysize=5cm  
\centerline{\epsfbox{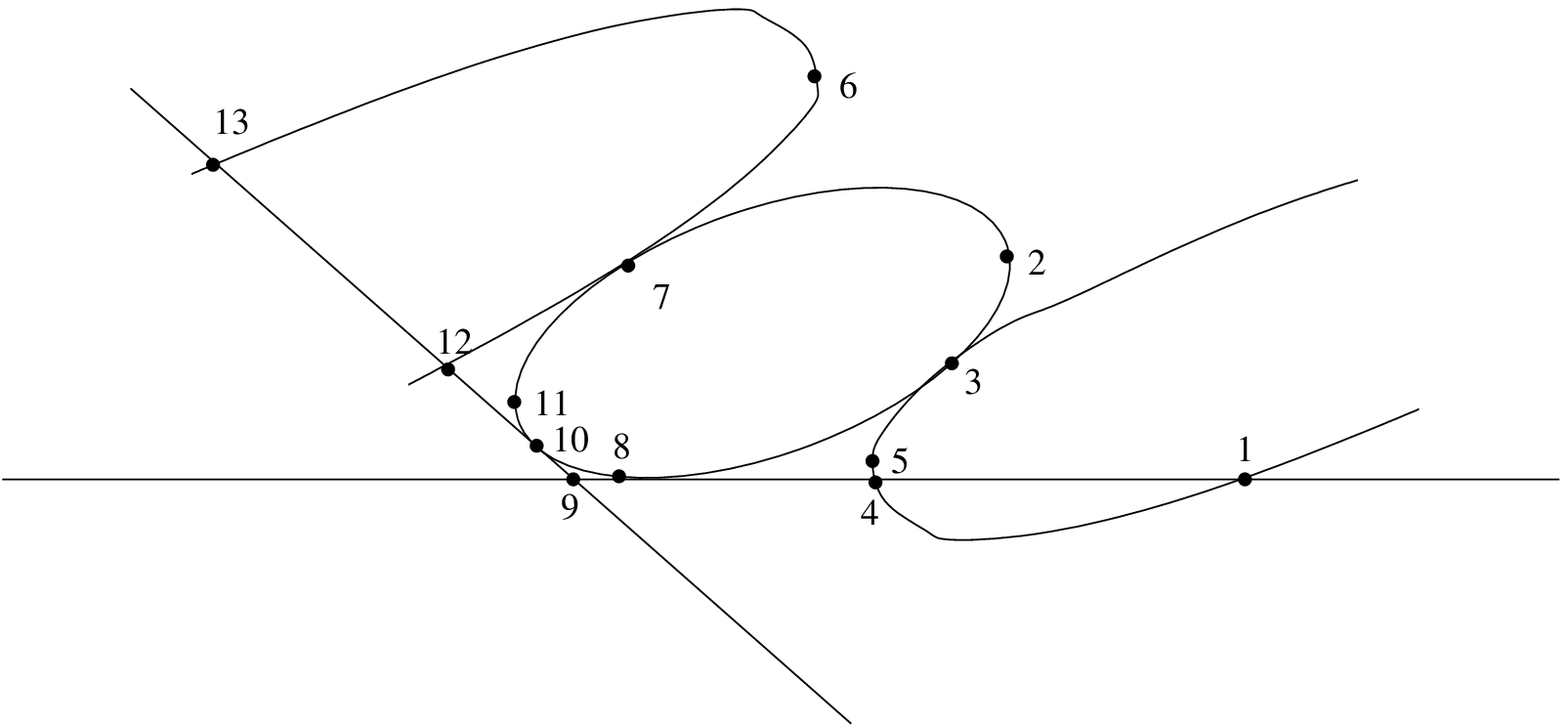}}  
\caption{The conic-line arrangement where both lines are tangent to the same conic}\label{thm4-fig}  
\end{figure}

\begin{thm}
Let $S$ be a curve in $\cpt$ composed of two tangent conics 
and two additional lines, 
which are both tangent to the same conic
(see Figure \ref{thm4-fig}), then:
$$\pcpt \cong \left\langle 
\begin{array}{c|c} 
x_1,x_2,x_3 & (x_2 x_3)^2=(x_3 x_2)^2, (x_1 x_3)^2=(x_3 x_1)^2,\\
            & [x_1, x_2]=[x_2, x_3 x_1 x_3^{-1}]=e 
\end{array}
\right\rangle$$  
\end{thm}

The proofs of the above theorems are similar to the ones in Section \ref{2conics_1line}.

\subsection{Arrangements with other singularities}

Apart from the two previous cases of arrangements with only simple tangency points, 
we have $16$ more possibilities to locate two lines into 
a configuration of two tangent conics.

In the following subsections, we  compute only the cases which turn 
out to have nonisomorphic fundamental groups of the complements. 
In the appendix, we list all the 16 possibilities with the corresponding 
fundamental groups of their complements.

\subsubsection{First case}

\begin{figure}[h]
\epsfysize=4cm  
\centerline{\epsfbox{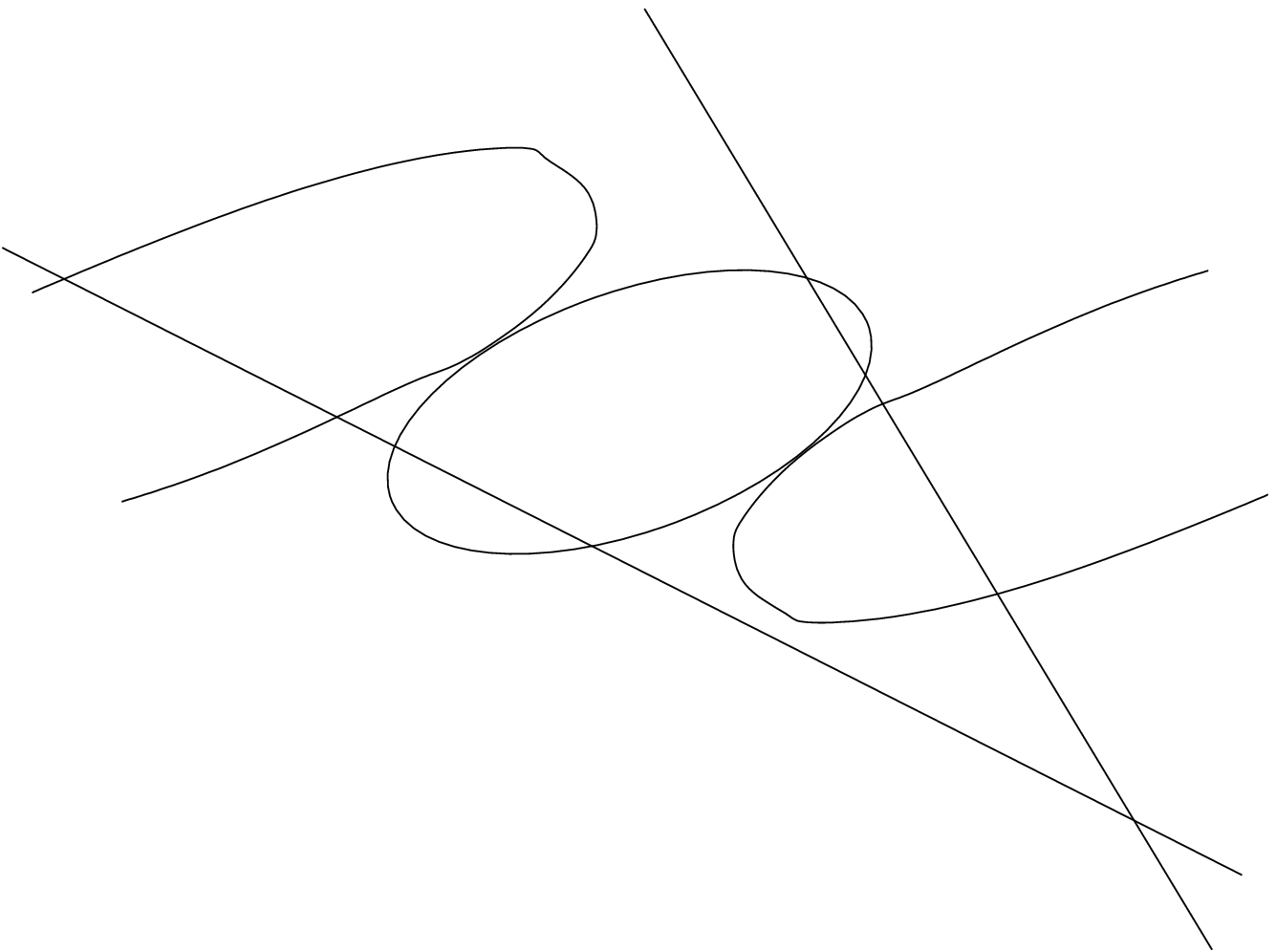}}
\caption{First case}\label{case2_16}  
\end{figure}

In this subsection, we  compute the fundamental group of the complement of 
the curve $S$ presented in Figure \ref{case2_16}.

Based on Proposition \ref{th1-2} and \cite{Ga}, 
the fundamental group of the complement of this case is  
$$\Z^2 \oplus \langle x_1,x_2 \ | \ (x_1 x_2)^2=(x_2 x_1)^2=e \rangle,$$
where $e$ is the identity element of the group.

\subsubsection{Second case}

\begin{figure}[h]
\epsfysize=4cm  
\centerline{\epsfbox{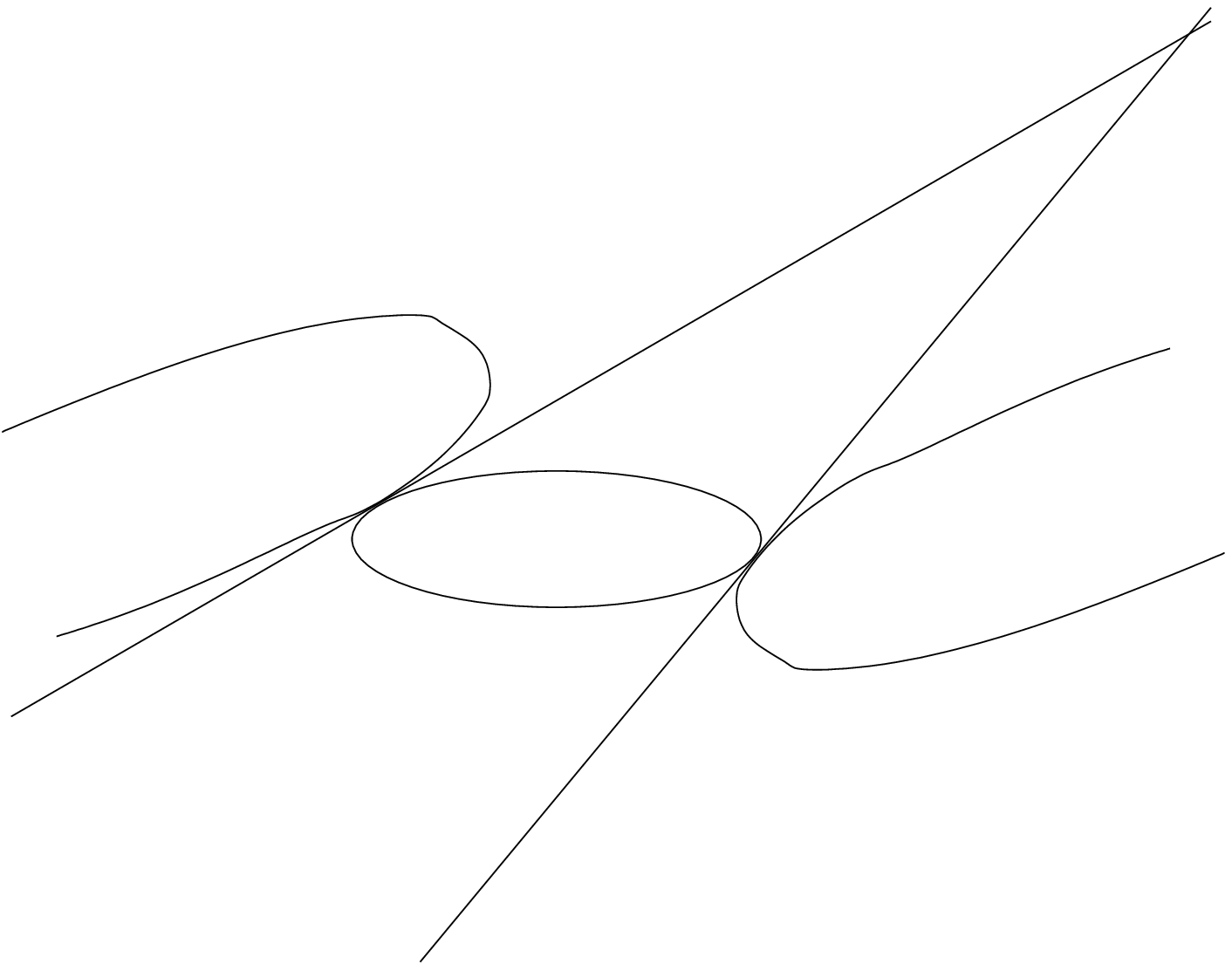}}  
\caption{Second case}\label{case2_1}  
\end{figure}

In this subsection, we  compute the fundamental group of the complement of 
the curve $S$ presented in Figure \ref{case2_1}.

Using braid monodromy techniques and the van Kampen theorem, 
we get the following presentation for the fundamental group of the complement:

\medskip

\noindent  
Generators: $\{ x_1,x_2,x_3,x_4,x_5,x_6 \}$. \\  
Relations:  
\begin{enumerate}    
\item $x_6 x_5 x_4 x_3 x_2 x_1 = e$ (projective relation)  
\item $x_3 x_5^{-1} x_6 x_5 = x_5^{-1} x_6 x_5 x_3$
\item $x_3^{-1} x_4 x_3 = x_5^{-1} x_6^{-1} x_5 x_6 x_5$
\item $(x_5 x_5^{-1} x_6 x_5 x_2)^2= (x_2 x_5 x_5^{-1} x_6 x_5)^2 = (x_5^{-1} x_6 x_5 x_2 x_5)^2$
\item $x_1 = x_6 x_5 x_2 x_5^{-1} x_6^{-1}$
\item $x_1 = x_4 x_3 x_2 x_3^{-1} x_4^{-1}$
\item $(x_4 x_3 x_2)^2 = (x_2 x_4 x_3)^2 = (x_3 x_2 x_4)^2$
\item $x_4 = x_5$
\end{enumerate}    

\medskip

In the first step, we use Relation (8) to substitute $x_4$ for $x_5$ 
everywhere. Moreover, since Relations (5) and (6) are both equal to $x_1$,
we can write the following equality:
$$x_6 x_4 x_2 x_4^{-1} x_6^{-1} =x_4 x_3 x_2 x_3^{-1} x_4^{-1}$$
Using Relation (5) again, we have that the projective relation 
gets the following form: 
$$x_6 x_4^2 x_3 x_2 x_6 x_4 x_2 x_4^{-1} x_6^{-1} =e$$
Hence, the generator $x_1$ is redundant.

So we get the following equivalent presentation:

\medskip
  
\noindent  
Generators: $\{ x_2,x_3,x_4,x_6 \}$. \\  
Relations:  
\begin{enumerate}    
\item $x_6 x_4^2 x_3 x_2 x_6 x_4 x_2 x_4^{-1} x_6^{-1} =e$  
\item $x_3 x_4^{-1} x_6 x_4 = x_4^{-1} x_6 x_4 x_3$
\item $x_3^{-1} x_4 x_3 = x_4^{-1} x_6^{-1} x_4 x_6 x_4$
\item $(x_6 x_4 x_2)^2= (x_2 x_6 x_4)^2 = (x_4^{-1} x_6 x_4 x_2 x_4)^2$
\item $x_6 x_4 x_2 x_4^{-1} x_6^{-1} = x_4 x_3 x_2 x_3^{-1} x_4^{-1}$
\item $(x_4 x_3 x_2)^2 = (x_2 x_4 x_3)^2 = (x_3 x_2 x_4)^2$
\end{enumerate}    

\medskip

From the first relation, we have that:
$x_4 x_3 x_2 x_6 x_4 x_2 =e$,   
and hence: $$x_3= x_4^{-1} x_2^{-1} x_4^{-1} x_6^{-1} x_2^{-1}$$
Now, we replace $x_3$ by $x_4^{-1} x_2^{-1} x_4^{-1} x_6^{-1} x_2^{-1}$, in any
place it appears.

Relation (2) gets the following form: 
$$x_4^{-1} x_2^{-1} x_4^{-1} x_6^{-1} x_2^{-1} x_4^{-1} x_6 x_4 = x_4^{-1} x_6 x_4 x_4^{-1} x_2^{-1} x_4^{-1} x_6^{-1} x_2^{-1}$$
By some reorderings, we get that:
$$(x_6 x_4 x_2)^2= (x_4 x_2 x_6)^2.$$

Relation (3) can now be written as:
$$x_2 x_6 x_4 x_2 x_4 x_4 x_4^{-1} x_2^{-1} x_4^{-1} x_6^{-1} x_2^{-1} = x_4^{-1} x_6^{-1} x_4 x_6 x_4$$
By some cancellations, we get that: $(x_6 x_4 x_2)^2 x_4 = x_4 (x_6 x_4 x_2)^2$.
By Relation (2), this is equal to: $(x_2 x_6 x_4)^2=(x_6 x_4 x_2)^2$, which 
is known (Relation (4)), and hence Relation (3) is redundant. 

Relation (5) can be written:
$$x_6 x_4 x_2 x_4^{-1} x_6^{-1} = x_4 x_4^{-1} x_2^{-1} x_4^{-1} x_6^{-1} x_2^{-1} x_2 x_2 x_6 x_4 x_2 x_4 x_4^{-1}$$
By some cancellations and reorderings, this is equal to:
$$(x_6 x_4 x_2)^2 = (x_2 x_6 x_4)^2$$
which is known (from Relation (2)), and hence this relation is redundant too.

By similar computations, Relation (6) can be written as:
$$(x_2 x_6 x_4)^2= (x_6 x_4 x_2)^2= (x_4 x_2 x_6)^2$$
which is known (Relations (2) and (4)), and hence Relation (6) is 
redundant too. 

The right part of Relation (4) can be written as:
$$(x_2 x_6 x_4)^2 = x_4^{-1} (x_6 x_4 x_2)^2 x_4$$
which is known by Relation (2), and therefore it is redundant too.

Hence, we get the following equivalent presentation: 
$$\langle  x_2,x_4,x_6 \  |\  (x_6 x_4 x_2)^2= (x_2 x_6 x_4)^2 = (x_4 x_2 x_6)^2 \rangle$$
as needed.

\subsubsection{Third case}
\begin{figure}[h]
\epsfysize=4cm  
\centerline{\epsfbox{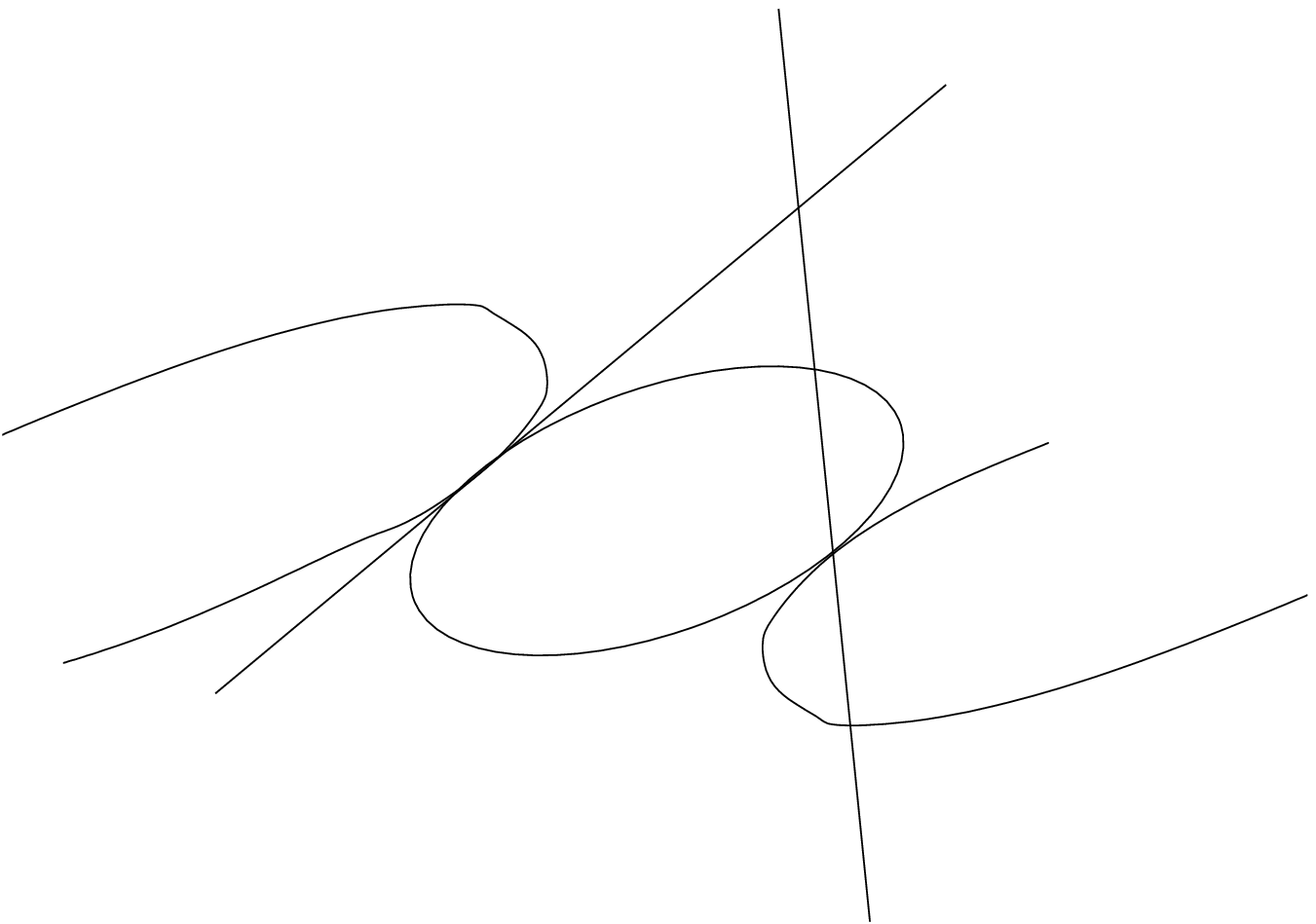}}  
\caption{Third case}\label{case2_2}  
\end{figure}

In this subsection, we  compute the fundamental group of the complement of 
the curve $S$ presented in Figure \ref{case2_2}.

Using braid monodromy techniques and the van Kampen theorem, 
we get the following presentation for the fundamental group of the complement:

\noindent  
Generators: $\{ x_1,x_2,x_3,x_4,x_5,x_6 \}$. \\  
Relations:  
\begin{enumerate}    
\item $x_6 x_5 x_4 x_3 x_2 x_1 = e$ (projective relation)  
\item $x_4 = x_5$ 
\item $x_1 x_2 = x_2 x_1$
\item $x_2 x_1 x_2^{-1} x_5 x_3 = x_5 x_3 x_2 x_1 x_2^{-1}$
\item $x_5 x_3 x_2 x_1 x_2^{-1} x_5 x_3 = x_3 x_5 x_3 x_2 x_1 x_2^{-1} x_5$ 
\item $x_1 x_4 = x_4 x_1$
\item $x_1 x_6 = x_6 x_1$
\item $x_2 = x_5 x_3 x_5^{-1}$
\item $x_2 = x_6 x_4 x_3 x_4^{-1} x_6^{-1}$ 
\item $(x_6 x_5 x_4 x_5^{-1} x_6^{-1} x_6 x_5 x_3 x_5^{-1})^2 = (x_5 x_3 x_5^{-1} x_6 x_5 x_4 x_5^{-1} 
x_6^{-1} x_6)^2 =\\ (x_6 x_5 x_3 x_5^{-1} x_6 x_5 x_4 x_5^{-1} x_6^{-1})^2$
\item $x_4 = x_6 x_5 x_6^{-1}$ 
\end{enumerate}    

Using Relations (2) and (8), we can cancel generators $x_2$ and $x_5$, and replace them 
by $x_4 x_3 x_4^{-1}$ and $x_4$ respectively. 
So we get the following simplified presentation:

\noindent  
Generators: $\{ x_1,x_3,x_4,x_6 \}$. \\  
Relations:  
\begin{enumerate}    
\item $x_6 x_4 x_4 x_3 x_4 x_3 x_4^{-1} x_1 = e$  
\item $x_1 x_4 x_3 x_4^{-1} = x_4 x_3 x_4^{-1} x_1$
\item $x_4 x_3 x_4^{-1} x_1 x_4 = x_4 x_3 x_4 x_3 x_4^{-1} x_1 x_4 x_3^{-1} x_4^{-1}$
\item $x_4 x_3 x_4 x_3 x_4^{-1} x_1 x_4 = x_3 x_4 x_3 x_4 x_3 x_4^{-1} x_1 x_4 x_3^{-1}$ 
\item $x_1 x_4 = x_4 x_1$
\item $x_1 x_6 = x_6 x_1$
\item $x_4 x_3 x_4^{-1} = x_6 x_4 x_3 x_4^{-1} x_6^{-1}$ 
\item $(x_6 x_4 x_6^{-1} x_6 x_4 x_3 x_4^{-1})^2 = (x_4 x_3 x_4^{-1} x_6 x_4 x_6^{-1} x_6)^2 = (x_6 x_4 x_3 x_4^{-1} x_6 x_4 x_6^{-1})^2$
\item $x_4 = x_6 x_4 x_6^{-1}$ 
\end{enumerate}    

By Relation (5), Relation (2) is reduced to $x_4 x_3 = x_3 x_1$. By this relation and 
Relation (5) again, Relation (3) can be simplified to the trivial relation, and hence it is 
redundant.
Moreover, Relation (4) gets the following form: $(x_4 x_3)^2 = (x_3 x_4)^2$.

By Relation (9), we can simplify Relations (7) and (8): they get the following forms respectively:
$x_3 x_6 = x_6 x_3$ and $(x_4 x_6 x_4 x_3 x_4^{-1})^2 = (x_4 x_3 x_6)^2 = (x_6 x_4 x_3)^2$.

Hence, we get the following simplified presentation:

\noindent  
Generators: $\{ x_1,x_3,x_4,x_6 \}$. \\  
Relations:  
\begin{enumerate}    
\item $x_6 x_4 x_4 x_3 x_4 x_3 x_4^{-1} x_1 = e$  
\item $x_1 x_3 = x_3 x_1$
\item $(x_4 x_3)^2 = (x_3 x_4)^2$ 
\item $x_1 x_4 = x_4 x_1$
\item $x_1 x_6 = x_6 x_1$
\item $x_3 x_6 = x_6 x_3$ 
\item $(x_4 x_6 x_4 x_3 x_4^{-1})^2 = (x_4 x_3 x_6)^2 = (x_6 x_4 x_3)^2$
\item $x_4 x_6 = x_6 x_4$ 
\end{enumerate}    

Since $x_6$ commutes with all the other generators, Relation (7) is reduced to 
$(x_4^2 x_3 x_4^{-1})^2 = (x_4 x_3)^2 = (x_4 x_3)^2$ 
The right equation is known (Relation (3)), so we have to simplify only the left equation.
By some simplifications, we get $x_4 x_3 x_4 x_3 =  x_3 x_4 x_3 x_4$ which is known too, 
and hence this relation is redundant.

Now, by Relation (1), we have  that $x_6 =  x_1^{-1}(x_4 x_3)^{-2}$. So we can replace $x_6$ 
by $x_1^{-1}(x_4 x_3)^{-2}$ in any place it appears.

Relation (5) gets the form:
$x_1  x_1^{-1}(x_4 x_3)^{-2}=  x_1^{-1}(x_4 x_3)^{-2} x_1$
Since $x_1$ commutes with $x_3$ and $x_4$, this relation becomes trivial.

Relations (6) and (8) get the following forms respectively:
$x_3 x_1^{-1}(x_4 x_3)^{-2} = x_1^{-1}(x_4 x_3)^{-2} x_3$ and
$x_4 x_1^{-1}(x_4 x_3)^{-2} = x_1^{-1}(x_4 x_3)^{-2} x_4$. 
Again, by using Relations (2), (3) and (4), these relations become trivial.

Hence, we get the following simplified presentation:

\noindent  
Generators: $\{ x_1,x_3,x_4 \}$. \\  
Relations:  
\begin{enumerate}    
\item $x_1 x_3 = x_3 x_1$
\item $(x_4 x_3)^2 = (x_3 x_4)^2$ 
\item $x_1 x_4 = x_4 x_1$
\end{enumerate}    
 
This group is obviously isomorphic to:
$$\langle x_1 \rangle \oplus \langle x_3,x_4 \ | \ (x_4 x_3)^2 = (x_3 x_4)^2 \rangle$$
as needed.

\subsubsection{Fourth case}

\begin{figure}[h]
\epsfysize=4cm  
\centerline{\epsfbox{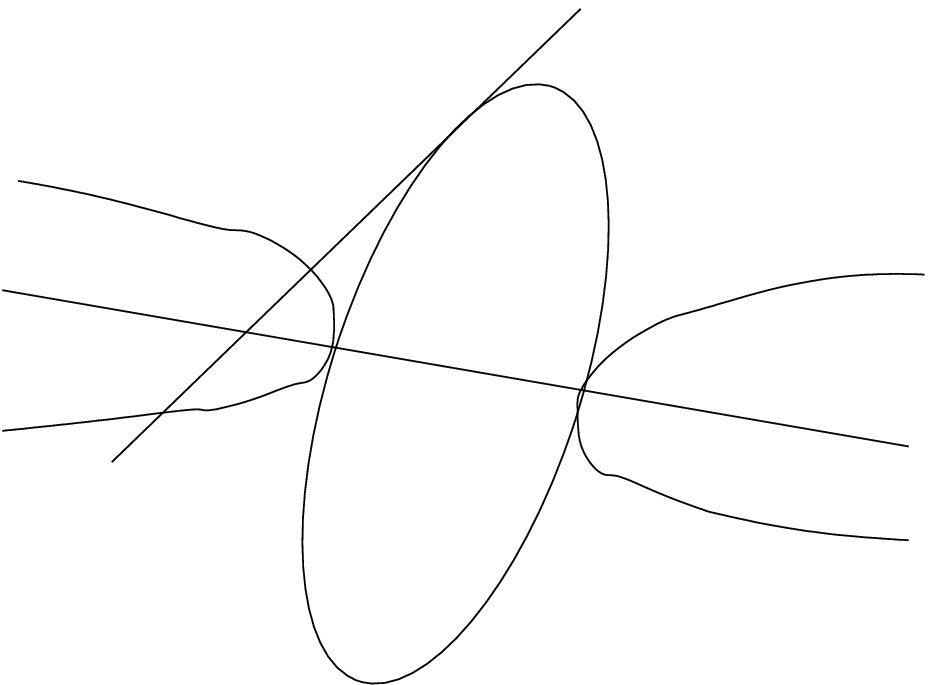}}  
\caption{Fourth case}\label{case2_7}  
\end{figure}

In this subsection, we  compute the fundamental group of the complement of 
the curve $S$ presented in Figure \ref{case2_7}.

Using braid monodromy techniques and the van Kampen theorem, 
we get the following presentation for the fundamental group of the complement:

\noindent  
Generators: $\{ x_1,x_2,x_3,x_4,x_5,x_6 \}$. \\  
Relations:  
\begin{enumerate}    
\item $x_6 x_5 x_4 x_3 x_2 x_1 = e$ (projective relation)  
\item $x_4 = x_5$
\item $x_1 = x_4 x_3 x_4^{-1}$
\item $x_4 x_3 x_2 = x_2 x_4 x_3$
\item $x_4 x_3 x_2 x_4 x_3 = x_3 x_2 x_4 x_3 x_4$
\item $(x_4 x_6)^2 = (x_6 x_4)^2$
\item $x_4 x_3 x_4^{-1} x_4 x_2 x_4 x_3 x_4^{-1} x_4 = x_4 x_2 x_4 x_3 x_4^{-1} x_4 x_4 x_3 x_4^{-1}$
\item $x_2 x_4 x_3 x_4^{-1} x_4 = x_4 x_3 x_4^{-1} x_4 x_2$
\item $x_4 x_3 x_4^{-1} = x_2^{-1} x_3^{-1} x_4^{-1} x_6 x_4^{-1} x_6^{-1} x_4 x_3 x_4^{-1} x_6 x_4 x_6^{-1} x_4 x_3 x_2$
\item $x_2^{-1} x_3^{-1} x_4^{-1} x_6 x_4 x_6^{-1} x_4 x_3 x_2 = x_4 x_3 x_4^{-1} x_4 x_4 x_3^{-1} x_4^{-1}$
\item $x_4 x_1 x_4^{-1} x_6 = x_6 x_4 x_1 x_4^{-1}$ 
\item $x_4 x_2 x_4^{-1} x_6 = x_6 x_4 x_2 x_4^{-1}$ 
\item $x_5 x_3 x_5^{-1} x_6 = x_6 x_5 x_3 x_5^{-1}$
\end{enumerate}    

Using Relations (2) and (3), we can cancel generators $x_1$ and $x_5$, and replace them 
by $x_4 x_3 x_4^{-1}$ and $x_4$ respectively. By a trivial simplification, Relation (8) 
is equal to Relation (4), so it is redundant.
So we get the following simplified presentation:

\noindent  
Generators: $\{ x_2,x_3,x_4,x_6 \}$. \\  
Relations:  
\begin{enumerate}    
\item $x_6 x_4^2 x_3 x_2 x_4 x_3 x_4^{-1} = e$   
\item $x_4 x_3 x_2 = x_2 x_4 x_3$
\item $x_4 x_3 x_2 x_4 x_3 = x_3 x_2 x_4 x_3 x_4$
\item $(x_4 x_6)^2 = (x_6 x_4)^2$
\item $x_3 x_2 x_4 x_3 = x_2 x_4 x_3 x_4 x_3 x_4^{-1}$
\item $x_4 x_3 x_4^{-1} = x_2^{-1} x_3^{-1} x_4^{-1} x_6 x_4^{-1} x_6^{-1} x_4 x_3 x_4^{-1} x_6 x_4 x_6^{-1} x_4 x_3 x_2$
\item $x_2^{-1} x_3^{-1} x_4^{-1} x_6 x_4 x_6^{-1} x_4 x_3 x_2 = x_4 x_3 x_4 x_3^{-1} x_4^{-1}$
\item $x_4^2 x_3 x_4^{-2} x_6 = x_6 x_4^2 x_3 x_4^{-2}$ 
\item $x_4 x_2 x_4^{-1} x_6 = x_6 x_4 x_2 x_4^{-1}$ 
\item $x_4 x_3 x_4^{-1} x_6 = x_6 x_4 x_3 x_4^{-1}$
\end{enumerate}    

From Relation (5) we have $x_3 x_2 x_4 x_3 x_4 = x_2 x_4 x_3 x_4 x_3$ which is equal to 
Relation (3) by Relation (2). Hence, Relation (5) is redundant.

Also, Relation (8) follows immediately from Relation (10), and hence it is redudant too.

Since the right side of Relation (1) is $e$, we can move $x_4 x_3 x_4^{-1}$ to the left, 
in order to get: $x_4 x_3 x_4^{-1} x_6 x_4^2 x_3 x_2 = e$. By Relation (10) and a cancellation,
we have $x_6 (x_4 x_3)^2 x_2 = e$.

Hence, we get the following presentation:  

\noindent  
Generators: $\{ x_2,x_3,x_4,x_6 \}$. \\  
Relations:  
\begin{enumerate}    
\item $x_6 (x_4 x_3)^2 x_2 = e$   
\item $x_4 x_3 x_2 = x_2 x_4 x_3$
\item $x_4 x_3 x_2 x_4 x_3 = x_3 x_2 x_4 x_3 x_4$
\item $(x_4 x_6)^2 = (x_6 x_4)^2$
\item $x_4 x_3 x_4^{-1} = x_2^{-1} x_3^{-1} x_4^{-1} x_6 x_4^{-1} x_6^{-1} x_4 x_3 x_4^{-1} x_6 x_4 x_6^{-1} x_4 x_3 x_2$
\item $x_2^{-1} x_3^{-1} x_4^{-1} x_6 x_4 x_6^{-1} x_4 x_3 x_2 = x_4 x_3 x_4 x_3^{-1} x_4^{-1}$
\item $x_4 x_2 x_4^{-1} x_6 = x_6 x_4 x_2 x_4^{-1}$ 
\item $x_4 x_3 x_4^{-1} x_6 = x_6 x_4 x_3 x_4^{-1}$
\end{enumerate}    

By the first relation $x_2 = (x_4 x_3)^{-2} x_6^{-1}$, so we can replace $x_2$ by $(x_4 x_3)^{-2} x_6$
in any place it appears. We start with Relation (7). We have:
$x_4 (x_4 x_3)^{-2} x_6^{-1} x_4^{-1} x_6 = x_6 x_4 (x_4 x_3)^{-2} x_6^{-1} x_4^{-1}$.  
By some simplifications, we have:
$x_6 x_4 x_6 x_4 x_3 x_4 x_3 = x_4 x_6 x_4 x_3 x_4 x_3 x_4^{-1} x_6 x_4$ 
By Relations (4) and (8), we have: $x_6 x_3 = x_3 x_6$.

Relation (2) gets the following form: 
$x_4 x_3 (x_4 x_3)^{-2} x_6^{-1} = (x_4 x_3)^{-2} x_6^{-1} x_4 x_3$. By some simplifications,
we have: $x_4 x_3 x_6 = x_6 x_4 x_3$. By the relation $x_6 x_3 = x_3 x_6$, we finally have 
that: $x_4 x_6 = x_6 x_4$. Therefore, Relation (4) is now redundant. Moreover, Relation (8) 
can be reduced to $x_3 x_6= x_6 x_3$, which is known, and hence Relation (8) is redundant too.

Since $x_3$ and $x_4$ commute with $x_6$, Relations (5) and (6) get the following forms 
respectively:
$x_4 x_3 x_4^{-1} = (x_6 (x_4 x_3)^2) x_3^{-1} x_4^{-1} x_3 x_4 x_3 ((x_4 x_3)^{-2} x_6^{-1})$ and 
$(x_6 (x_4 x_3)^2) x_3^{-1} x_4 x_3 ((x_4 x_3)^{-2} x_6^{-1}) = x_4 x_3 x_4 x_3^{-1} x_4^{-1}$.
If we continue the simplifications, we get that both relations become trivial, 
and hence they are redundant.

Finally, we simplify Relation (3). By Relation (1), it gets the following form: 
$x_4 x_3 ((x_4 x_3)^{-2} x_6^{-1}) x_4 x_3 = x_3 ((x_4 x_3)^{-2} x_6^{-1}) x_4 x_3 x_4$. 
By some simplifications, we get that this relation becomes trivial, and hence it is redundant too.

Therefore, we get the following simplified presentation: 

\noindent  
Generators: $\{ x_3,x_4,x_6 \}$. \\  
Relations:  
\begin{enumerate}    
\item $x_3 x_6 = x_6 x_3$   
\item $x_4 x_6 = x_6 x_4$
\end{enumerate}    

This group is obviously the group:
$$\Z \oplus \F_2$$
as needed. 

\section{The fundamental groups are big}\label{sec_big}
A group is called {\it big} if it contains a subgroup which is free 
(generated by two or more generators). 
In this section we  show that the fundamental groups which appear in this paper are big. 

We start with a simple observation. If $G$ is a group which has a big quotient, then 
$G$ is big itself: Let $N$ be a normal subgroup of $G$ such that $G/N$ is big. 
Let $aN,bN$ be the two generators of the free subgroup of $G/N$, then the group generated 
by $a,b$ in $G$ is big. Otherwise, there is a relation $w(a,b)=1$ in $G$, and therefore
we have a corresponding relation $w(aN,bN)=N$, which is a contradiction to the freeness of
the subgroup $\langle aN,bN \rangle \leq G/N$.   

\medskip

Using the observation, we show the following result:

\begin{prop}
The group
$$G \cong \langle a,b\ | \ (ab)^2 = (ba)^2 =e \rangle$$
is big.
\end{prop}

\begin{proof}
First, we change the presentation: let us take new generators $x=ab, y=b$, then the relation 
$(ab)^2 = (ba)^2 =e$ becomes: $x^2= (yxy^{-1})^2=e$, which is equal to: $x^2= y x^2 y^{-1}=e$.
Hence, we have:
$$G \cong \langle x,y\ |\ x^2= y x^2 y^{-1}=e \rangle$$
Now, the quotient of $G$ by the subgroup generated by $y^3$ is 
$$G/\langle\langle y^3 \rangle\rangle \cong \Z/2 * \Z/3$$
where $*$ is the free product. 

Since $\Z/2 * \Z/3$ is known to be big, by the observation we have that $G$ is big as needed.
\end{proof}

Now we can prove Corollary \ref{bigness}. Let $A$ be a conic-line arrangement with 
two tangent conics and up to two additional lines. 
Since the two tangent conics is a sub-arrangement of $A$, 
it is known that the fundamental group of the 
two tangent conics is a quotient of 
the fundamental group of $A$ (by sending  the generators 
which correspond to the additional lines to $e$). 
Since, we have shown in the previous proposition that the fundamental group of the 
two tangent conics is big, 
then we have that  the fundamental group of $A$ is big too 
(by using again the observation). Hence, we are done.

\section{Conjectures concerning the case of arrangements with only simple tangency points}\label{general} 
Motivated by the presentations we achieved in the cases of two tangent conics,
with up to two tangent lines with only simple tangnecy points, 
we want to propose some conjectures about 
the simplified presentation of the fundamental group of two tangent conics, 
with an arbitrary number of tangented lines, 
which are tangent in simple tangency points.

\begin{conj}
Let $S$ be a conic-line arrangement with two tangent conics, 
and $n$ lines which are tangent
to these conics. Assume that the tangency points of $S$ are simple. 
Let $C_1$ be the conic with the maximal number $m$ of lines which are 
tangent to it.
Then the following holds in the simplified presentation of $\pi_1(\C\P^2-S)$:
\begin{enumerate}
\item We have $m$ relations of the type
$$(a x_i)^2 = (x_i a)^2,$$ 
where $a$ is the generator corresponding to the conic $C_1$ and $x_i$ are 
some generators of the group. 
\item Let $y$ and $z$ be two generators of $\pi_1(\C\P^2-S)$ which are different from $a$. 
We have 
$$y aza^{-1} = aza^{-1} y$$
\end{enumerate}
\end{conj}

\begin{rem}
In some of the computed presentations, we do not have all the relations of the type mentioned
in the second part of the conjecture. This happens since some of these relations are 
simplified to commutative relations between generators by other relations.    
\end{rem}

For a specific case, we should get a decomposition:
\begin{conj}
Let $S$ be a conic-line arrangement with two tangent conics 
and $n$ lines, where $n-1$ 
lines are tangent to one conic, and one line $L$ is tangent to the other one.
Assume that all the tangency points are simple.
Let $x$ be the generator that corresponds to $L$.
Then:
$$\pi_1(\C\P^2-S) = \langle x \rangle \oplus \pi_1(\C\P^2-(S-L)).$$
\end{conj}

\section*{Appendix}
  
In the appendix, we list all the possibilities for two tangency points and 
two additional lines, and the corresponding fundmanetal groups of the complement
(we assume that the line at infinity is transversal to our curve, i.e. there are no
singularities on the line at infinity).

\medskip
\epsfysize=3cm  
\centerline{\epsfbox{case2_1.eps}}
$$\pi_1(\C\P^2-S) \cong \langle x,y,z | (zyx)^2=(yxz)^2=(xzy)^2 \rangle$$

\medskip
\epsfysize=10cm 
\epsfxsize=9cm 
\centerline{\epsfbox{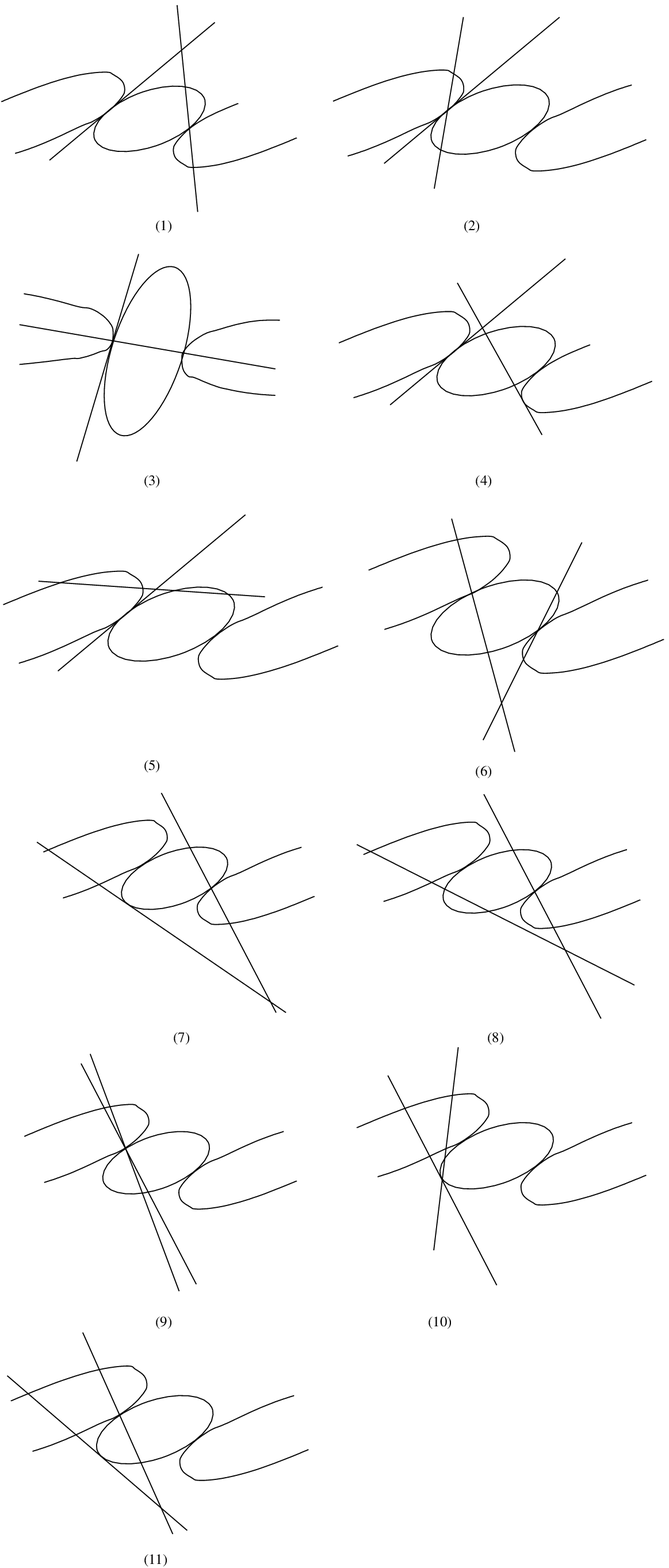}}
$$\pi_1(\C\P^2-S) \cong \langle x \rangle \oplus \langle y,z | (yz)^2=(zy)^2 \rangle$$

\medskip
\epsfysize=2.5cm
\epsfxsize=9cm 
\centerline{\epsfbox{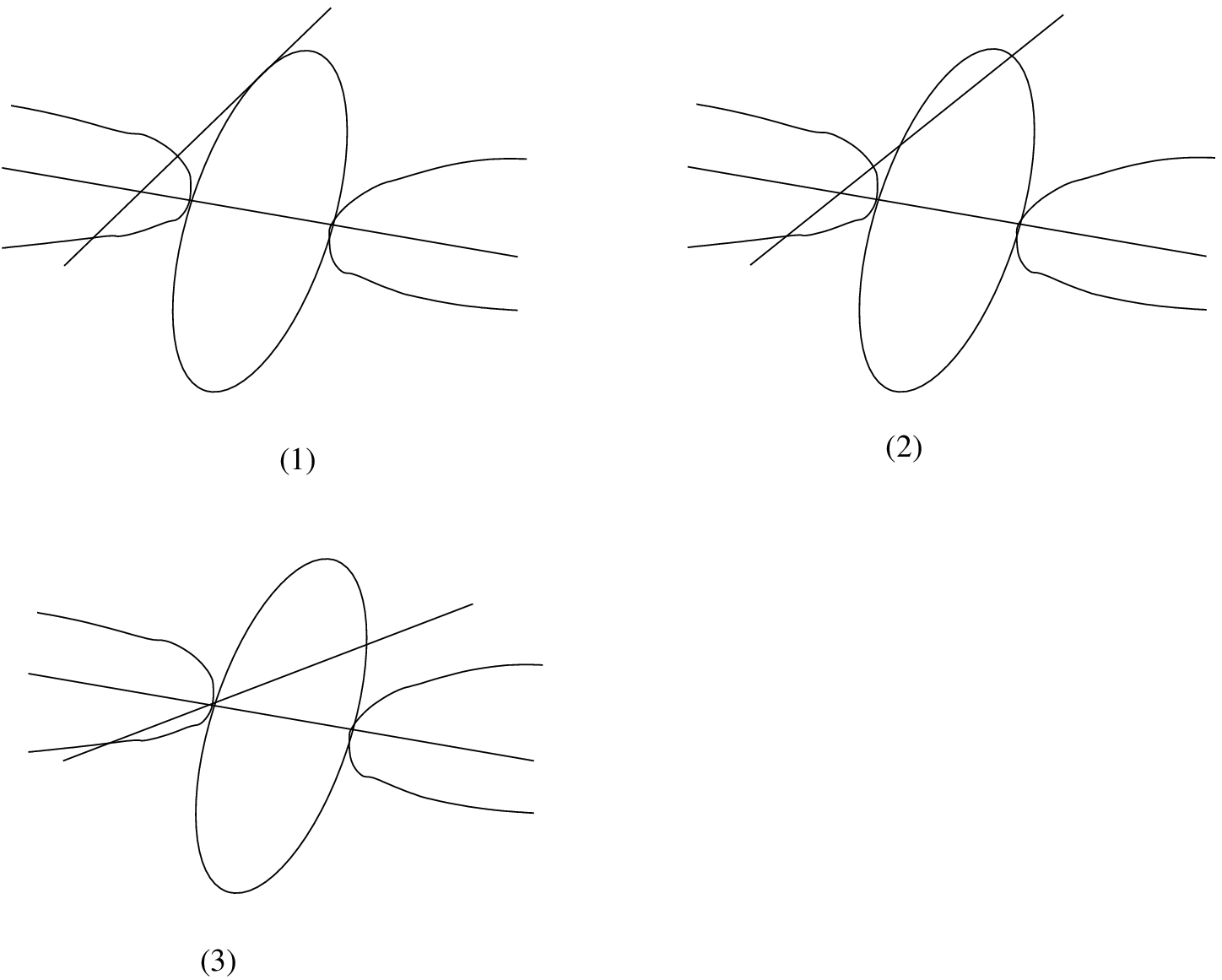}}  
$$\pi_1(\C\P^2-S) \cong \Z \oplus \F_2$$

\medskip
\epsfysize=1.5cm  
\centerline{\epsfbox{case2_16.eps}}
$$\pi_1(\C\P^2-S) \cong \Z ^2 \oplus  \langle y,z | (xy)^2=(yx)^2 =e \rangle$$

\section*{Acknowledgments} 

We wish to thank Alex Degtyarev for fruitful discussions.

The first and second authors wish to thank the Institute of Mathematics, 
Erlangen (Germany), Institut Fourier, Grenoble (France), Hershel Farkas, 
Ron Livne and Einstein Institute of Mathematics for hosting their stays.

\end{document}